\RequirePackage{etoolbox}
\csdef{input@path}%
{%
 {sty/},
 {img/},
}
\documentclass[a4paper,p04,,usenames,dvipsnames]{article}
\usepackage[numbers,square]{natbib} 
\usepackage{graphicx}
\usepackage{subfig}
\usepackage{acro}
\usepackage{geometry}
\usepackage{comment}
\usepackage{xspace}
\usepackage{amsmath,amssymb,bm}
\usepackage{longtable}
\usepackage{overpic}
\usepackage{tikz}
\usepackage{pict2e}

\usepackage[
]{caption}

\usepackage{hyperref}
\hypersetup{%
	pdfpagemode={UseOutlines},
	bookmarksopen,
	pdfstartview={FitH},
	colorlinks,
	linkcolor={cyan},
	citecolor={cyan},
	urlcolor={cyan}
}

\definecolor{emerald}{HTML}{50C878}
\definecolor{bronze}{HTML}{CD7F32}
\definecolor{brickred}{HTML}{CB4154}

\renewcommand*\thesection{\arabic{section}}
\renewcommand\bibname{REFERENCES}

\usepackage{algorithm}
\usepackage{algpseudocode}

{\sffamily\large\itshape\bfseries}

\def\eqref#1{\textcolor{cyan}{(\ref{#1})}}

\begin{document}

\newcommand{\openfoam}{Open\nolinebreak\hspace{-.2em}\nolinebreak\hspace{-.2em}FOAM\textsuperscript{\textregistered}\xspace}








\title{On the accuracy and efficiency of
reduced order models: towards real-world applications}
\author{Pierfrancesco Siena, Pasquale Claudio Africa, Michele Girfoglio, Gianluigi Rozza}
\date{}
\maketitle

\begin{abstract}
This chapter provides an extended overview about Reduced Order Models (ROMs), with a focus on their features in terms of efficiency and accuracy. In particular, the aim is to browse the more common ROM frameworks, considering both intrusive and data-driven approaches. We present the validation of such techniques against several test cases. The first one is an academic benchmark, the thermal block problem, where a Poisson equation is considered. Here a classic intrusive ROM framework based on a Galerkin projection scheme is employed. The second and third test cases come from real-world applications, the one related to the investigation of the blood flow patterns in a patient specific coronary arteries configuration where the Navier Stokes equations are addressed and the other one  
concerning the granulation process within pharmaceutical industry where a fluid-particle system is considered.  
Here we employ two data-driven ROM approaches showing a very relevant trade-off between accuracy and efficiency. 
In the last part of the contribution, two novel technological platforms, ARGOS and ATLAS, are presented. They are designed to provide a user-friendly access to data-driven models for real-time predictions for complex biomedical and industrial problems. 
\end{abstract}

\textbf{keywords:} Reduced Order Models,  Data-driven techiques, Real-time computing 

\section{Introduction}\label{sec:intro}
In the field of computational modeling, the trade-off between efficiency and accuracy is a central aspect for researchers and specialists. High-fidelity discretization methods (e.g.,
finite element, finite volume, or spectral element methods), hereinafter referred to as Full Order Models (FOMs) are able to offer insights into complex phenomena, however the cost of computational resources and time is prohibitive for real-world applications, especially when many physical or geometrical parameters are involved. 

Conversely, Reduced Order Models (ROMs) \cite{porsching1985estimation,ito1998reduced,lassila2013generalized,gunzburger2007reduced,rozza2008reduced,aromabook,hesthaven2022reduced,siena2023introduction} provide a framework to speed up FOMs but on the other hand they may raise concerns regarding the accuracy of the predicted solution. 
The capability of ROMs to achieve real-time evaluations derives from the observation that the parametric dependencies of the problem at hand, in many cases, has an intrinsic dimensionality significantly lower than the number of degrees of freedom associated with the FOM. 

We focus on the Reduced Basis (RB) method: 
the goal is to identify a small number of basis functions that can be linearly combined with respect to proper weights (the so called reduced coefficients) to accurately reconstruct the solution 
at a lower computational cost. 
Then a fully decoupled \emph{offline}-\emph{online} paradigm \cite{benner2020model,rozza20201} is introduced. In the \emph{offline} stage the FOM is solved for selected values of the parameters involved in the problem, so its computational cost is in general expensive. Anyway, this does not represent an issue because the collection of high-fidelity snapshots is performed only once. Then, a reduced basis (in principle of small size) is computed from the database of FOM solutions. 
On the other hand, in the \emph{online} stage, the prediction of the reduced coefficients is carried out usually either through interpolation/regression techniques or by solving a system of algebraic equations or ordinary differential equations (for example in the case of time dependent problems) obtained by projecting the FOM onto the reduced space spanned by the basis functions. The online stage enables the parameter space to be explored with a significantly reduced cost.

\textcolor{black}{Concerning the reliability of ROMs, approaches such as residual-based error estimation, a posteriori estimation, and uncertainty quantification can provide insights into errors in ROMs \cite{quarteroni2014reduced,quarteroni2011certified,feng2024posteriori,nurtaj2020adaptive}. These methods analyze the difference between the computed solution and the FOM solution, involve solving additional equations or systems to estimate the error or deal with uncertainties arising from various sources such as model parameters, boundary conditions, or numerical discretization.  However, several open problems still persist. These include the development of efficient techniques for grid adaptation, the integration of uncertainty into ROMs, especially in multi-physics or coupled scenarios, and the enhancment of error estimation capabilities to handle the inherent complexity of ROMs effectively.}

This chapter sets out to explore the balance between accuracy and efficiency, honing in on why reduced order models might have a real-world impact. Therefore, the interest is not just on how fast ROMs are, but also on how accurate they can be. 
By delving into the accuracy of ROMs beyond the computational features, we aim to test their reliability and robustness in real-world applications.


The rest of this chapter is structured as follows. Sec. \ref{sec:rom} is dedicated to the mathematical description of intrusive and data-driven ROMs, underlying advantages and disadvantages of both approaches. Then, an intrusive ROM framework is applied to an academic benchmark case while some data-driven ROM approaches for real world cases arising from biomedical and industrial engineering applications. In particular, an hemodynamics study is shown in Sec. \ref{sec:coronary} for coronary arteries and an industrial benchmark is proposed in Sec. \ref{sec:industry}. Finally, we present our two online computational platforms for real-time simulations, namely Argos and Atlas, in Sec. \ref{sec:argos_atlas}.

\section{Reduced order models}\label{sec:rom}

In this section, we provide a brief mathematical description about ROMs. Let us analyze in the following the two main categories of ROMs \cite{aromabook,benner2020model,benner2020model2,degruyter1,rozza20201}: intrusive ROMs, where the original governing equations play a fundamental role (see, e.g. \cite{stabile2017pod,stabile2018finite,girfoglio2021pod,girfoglio2021pressure,pichi2022driving,pichi2019reduced,pichi2020reduced}), and data-driven ROMs, which are based only on the data (see e.g. \cite{siena2023fast,Siena2022,balzotti2023reduced,zainib2022chapter,balzottidata2022,papapicco2021neural,hijazi2020data,shah2021finite,coscia2023physics,pichi2023artificial}). In particular, concerning intrusive ROMs, we refer to the popular Galerkin projection methods, where the governing equations are projected onto the reduced space by providing, in the most general framework, a differential algebraic system of gives whose unknows are the reduced coefficients \cite{rowley2004model,iollo2000stability}.   

Let us introduce the computational domain $\Omega \subset \mathbf{R}^d$ where $d=1,2,3$, and  its boundary $\partial \Omega$. The regions of $\partial \Omega$ where Dirichlet boundary conditions are enforced are indicated with $\Gamma^D_i$, where $1\leq i \leq d_v$. Let be $w:\Omega \mapsto \mathbf{R}^{d_v}$ a vector field ($d_v = d$) or a scalar field ($d_v = 1$).



\subsection{Galerkin projection ROMs}
\label{sec:intrusive-rom} 

For sake of simplicity we consider a stationary problem. Moreover, we refer to a finite element formulation of the problem at hand. Then the following scalar space is introduced:
\begin{equation}
    \mathbf{V}_i(\Omega)=\{ v\in H^1(\Omega) : v_{\mid \Gamma^D_i}=0 \}, \qquad  1\leq i \leq d_v,
\end{equation}
where $H^1_0(\Omega) \subset \mathbf{V}_i(\Omega) \subset H^1(\Omega)$. It is well known that if only Dirichlet boundary conditions are imposed, i.e. $ \Gamma^D_i = \partial \Omega $, the scalar space $\mathbf{V}_i(\Omega)$ coincides with the popular Hilbert space $H^1_0(\Omega)$. From here on out, we mention $\mathbf{V}_i(\Omega)$ as $\mathbf{V}_i$, for a smoother explanation.

Let us consider the product space $\mathbf{V}= \mathbf{V}_1 \times \dots \times \mathbf{V}_{d_v}$ and a general element $w = (w_1, \dots, w_{d_v}) \in \mathbf{V}$. The inner product identified by $(w,v)_{\mathbf{V}}, \forall w,v \in \mathbf{V}$ is introduced and the resulted norm is $\lVert w \rVert_{\mathbf{V}}=(w,w)_{\mathbf{V}}^{1/2}, \forall w \in \mathbf{V}$. Note that $\mathbf{V}$ is an Hilbert space.

In order to analyze parametric cases, a closed parameter domain $\mathbf{P} \subset \mathbf{R}^p$ is introduced. A general element of $\mathbf{P}$ is identified by $\mu = (\mu_1, \dots, \mu_p)$, so the following parametric field variable for a parameter value is provided  $u(\mu)=(u_1 (\mu),\dots, u_{d_v}(\mu)):\mathbf{P}\mapsto\mathbf{V}$.

Given $\mu \in \mathbf{P} \subset \mathbf{R}^p$, the abstract formulation of a stationary problem evaluates $s(\mu) = l(u(\mu);\mu)$ where $u(\mu)\in \mathbf{V}$ satisfies:
\begin{equation}
      a(u(\mu),v;\mu) = f(v;\mu) \qquad \forall v \in \mathbf{V}.
      \label{varitional-problem}
\end{equation}
In eq. \eqref{varitional-problem} $a:\mathbf{V}\times\mathbf{V}\times\mathbf{P}\mapsto\mathbf{R}$ is a bilinear form with respect to $\mathbf{V} \times \mathbf{V}$, $f:\mathbf{V}\times\mathbf{P}\mapsto\mathbf{R}$ and $l:\mathbf{V}\times\mathbf{P}\mapsto\mathbf{R}$ are linear forms with respect to $\mathbf{V}$ and $s:\mathbf{P}\mapsto\mathbf{R}$ is the scalar output of the model. The problem \eqref{varitional-problem} is our FOM in continuous form.

We shall presume for most of this analysis that the compliant hypotheses is verified. In particular we assume:
\begin{itemize}
    \item $f(\cdot;\mu) = l(\cdot;\mu) \quad \forall \mu \in \mathbf{P}$,
    \item $a(\cdot,\cdot;\mu)$ symmetric $\forall \mu \in \mathbf{P}$.
\end{itemize}
It greatly simplifies the following presentation while still showing most of the significant RB principles. 
We shall further assume the well-posedness of the problem \eqref{varitional-problem}. In particular, the Lax-Milgram theorem \cite{quarteroni2008numerical} states that if 
\begin{itemize}
    \item the bilinear form $a(\cdot,\cdot;\mu)$ is coercive and continuous on $\mathbf{V} \times \mathbf{V}$,
    \item the linear form $f(\cdot;\mu)$ is continuous on $\mathbf{V}$,
\end{itemize}
the problem \eqref{varitional-problem} admits a unique solution.
We shall adopt one last assumption, fundamental for the performance of the ROM framework in terms of efficiency: the affine decomposition. Moreover, this hypotheses allow to introduce the \emph{offline}-\emph{online} paradigm discussed in Sec. \ref{sec:intro}. 
The bilinear and linear forms above introduced are affine if they can be written as:
\begin{equation}
    a(u,v;\mu)=\sum_{q=1}^{Q_a} \theta^q_a(\mu)a_q(u,v), \quad
    f(v;\mu) = \sum_{q=1}^{Q_f}\theta^q_f(\mu) f_q(v), \quad
    l(v;\mu)=\sum_{q=1}^{Q_l} \theta_l^q(\mu)l_q(v),
    \label{affine}
\end{equation}
where $a_q$, $f_q$ and $l_q$ are independent on the parameter $\mu \in \mathbf{P}$ while $\theta^q_a$, $\theta^q_f$, $\theta^q_l$ are scalar quantities depending only on the parameter values $\mu \in \mathbf{P}$. \textcolor{black}{The reader is referred to Sec. \ref{sec:termal} for an illustrative example.} Notice that when the problem does not provide an affine decomposition, it can be approximated using other methods, as the empirical interpolation procedure \cite{hesthaven2016certified} (\textcolor{black}{for a deeper analysis of this procedure, readers are directed to the chapter by Nguyen \cite{nguyen2024model}}).

In order to discretize the problem \eqref{varitional-problem}, let introduce a finite dimensional subset $\mathbf{V}_{\delta} \subset \mathbf{V}$ with size $\text{dim}(\mathbf{V}_{\delta}) = N_\delta$. The finite dimensional form of equation \eqref{varitional-problem} evaluates $s_{\delta}(\mu) = l(u_{\delta}(\mu);\mu)$, where $u_{\delta}(\mu)$ satisfies
\begin{equation}
      a(u_{\delta}(\mu),v_{\delta};\mu) = f(v_{\delta};\mu) \qquad \forall v_{\delta} \in \mathbf{V}_{\delta}.
      \label{varitional-problem-discr}
\end{equation}
The problem \eqref{varitional-problem-discr} is our FOM in discrete form.
Thanks to the continuity, the coercivity of $a(\cdot,\cdot;\mu)$,  and the Galerkin orthogonality, the Cea's lemma \cite{monk2003finite} holds:
\begin{equation}
    \lVert u(\mu)- u_{\delta}(\mu) \rVert_{\mathbf{V}} \le \bigg( 1 + C(\mu)
    \bigg) \inf_{v_{\delta}\in \mathbf{V}_{\delta}} \lVert u(\mu)-v_{\delta} \rVert_{\mathbf{V}}, \quad \forall v_{\delta} \in \mathbf{V}_{\delta}, \label{eq:Cea1}
\end{equation}
where $C(\mu)$ depends on the coercivity and continuity constant of the form $a(\cdot,\cdot;\mu)$. As a result, the best approximation error of $u(\mu)$ in the approximation space $\mathbf{V}_{\delta}$ is strongly connected to the approximation error $\lVert u(\mu)- u_{\delta}(\mu) \rVert_{\mathbf{V}}$ through the constant $C(\mu)$.

Let introduce the solution manifold of the problem \eqref{varitional-problem}:
\begin{equation}
    \mathcal{M} = \{u(\mu) : \mu \in \mathbf{P} \} \subset \mathbf{V},
\end{equation}
collecting all the exact solutions $u(\mu)\in\mathbf{V}$ varying the parameter $\mu$.
Similarly, at a discrete level, the solution manifold of the problem \eqref{varitional-problem-discr} is:
\begin{equation}
    \mathcal{M}_{\delta}=\{ u_{\delta}(\mu) : \mu \in \mathbf{P}\} \subset \mathbf{V}_\delta.
\end{equation}

High fidelity solutions are typically only computed when the number of parameters is limited because of the typical high number of degrees of freedom $N_{\delta}$. 
The goal is to find a small number of basis functions whose linear combination exactly represents the numerical solution $u_{\delta}(\mu)$. Let $\{ \xi \}_{i=1}^{N}$ be an $N$-dimensional set of reduced basis (the most important techniques for its calculation are explained in Section \ref{sec:pod}). Then, the space derived is:
\begin{equation}
    \mathbf{V}_{\text{rb}}=\text{span}\{ \xi_1,\dots,\xi_N \} \subset \mathbf{V}_{\delta},
\end{equation}
where we suppose $N \ll N_{\delta}$. Therefore, the reduced form of equation \eqref{varitional-problem-discr} evaluates $s_{\text{rb}}(\mu) = f(u_{\text{rb}}(\mu);\mu)$, where $u_{\text{rb}}(\mu)$ satisfies:
  \begin{equation}
      a(u_{\text{rb}}(\mu),v_{\text{rb}};\mu) = f(v_{\text{rb}};\mu) \qquad \forall v_{\text{rb}} \in \mathbf{V}_{\text{rb}}.
      \label{varitional-problem-rb}
  \end{equation}
The reduced solution can be expressed as 
  \begin{equation}
u_{\text{rb}}(\mu)=\sum_{i=1}^{N}(u_{\text{rb}}^{\mu})_i \xi_i,
      \label{varitional-problem-rb2}
  \end{equation}
 where $(u_{\text{rb}}^{\mu})_i$ are the coefficients of the reduced basis approximation. Note that the reduced coefficients depend on the parameter $\mu$, whereas the reduced basis is independent on it. 

Now we  analyze the accuracy of the reduced solution. Using the triangle inequality, we have 
\begin{equation}
    \lVert u(\mu)- u_{\text{rb}}(\mu) \rVert_{\mathbf{V}} \le \lVert  u(\mu) - u_{\delta}(\mu) \rVert_{\mathbf{V}} + \lVert u_{\delta}(\mu) - u_{\text{rb}}(\mu)  \rVert_{\mathbf{V}}.
\end{equation}
A measure of the distance between $\mathcal{M}_{\delta}$ and $\mathbf{V}_{\text{rb}}$ was already introduced in literature, namely the Kolmogorov $N$-width \cite{pinkus2012n}, defined as:
\begin{equation}
    d_N(\mathcal{M}_{\delta}) = \inf_{\mathbf{V}_{\text{rb}}} \mathcal{E}(\mathcal{M}_{\delta},\mathbf{V}_{\text{rb}}) = \inf_{\mathbf{V}_{\text{rb}}} \sup_{u_{\delta}\in\mathcal{M}_{\delta}}\inf_{v_{\text{rb}}\in\mathbf{V}_{\text{rb}}} \lVert u_{\delta}-v_{\text{rb}} \rVert_{\mathbf{V}}.
\end{equation}
\textcolor{black}{This quantity is a measure used to understand how well the details of the system are captured. For example, in turbulent situations, essentially the flow is splitted into a defined number ($N$) of finely divided components to ensure accurate depiction. When the Kolmogorov $N$-width increases, a simulation gains more detail, but at the expense of greater computational demands. Achieving the right balance between accuracy and computational efficiency is key in modeling approach.} For elliptic problems, such as diffusion equations, the quantity $d_N(\mathcal{M}_{\delta})$ decreases quickly as $N$ increases. So, in these cases, a small number of basis functions can approximate accurately the set of high fidelity solutions $\mathcal{M}_{\delta}$ \cite{pinkus2012n,chen2012certified}. 
The Cea's lemma (see eq. \eqref{eq:Cea1}) for $\mathbf{V}_{\text{rb}}$ reads:
\begin{equation}
    \lVert u(\mu)- u_{\text{rb}}(\mu) \rVert_{\mathbf{V}} \le \bigg( 1 + 
    C(\mu)\bigg) \inf_{v_{\text{rb}}\in \mathbf{V}_{\text{rb}}} \lVert u(\mu)-v_{\text{rb}} \rVert_{\mathbf{V}}, \quad \forall v_{\text{rb}} \in \mathbf{V}_{\text{rb}}.
    \label{cea}
\end{equation}
Thus we can conclude that the ROM accuracy depends on the problem at hand.

However, for simple problems, it is possible to control the error during the online phase by means of some a posteriori estimations: 
\begin{equation}
    a(u_{\delta}(\mu)-u_{\text{rb}}(\mu),v_{\delta};\mu) = r(v_{\delta};\mu) = (\hat{r}_{\delta}(\mu),v_{\delta})_{\mathbf{V}}, \quad \forall v_{\delta} \in \mathbf{V}_{\delta},
    \label{error-eq1}
\end{equation}
where $r(\cdot,\mu)$ is the residual in the dual space $\mathbf{V}_{\delta}^*$  and $\hat{r}_{\delta}(\mu) \in \mathbf{V}_{\delta}$ is its Riesz representation \cite{hesthaven2016certified}.
Therefore it holds:
\begin{equation}
    \lVert \hat{r}_{\delta}(\mu) \rVert_{\mathbf{V}} = \lVert r(\cdot;\mu) \rVert_{\mathbf{V}_{\delta}^*} =  \sup_{v_{\delta}\in\mathbf{V}_{\delta}}\frac{r(v_{\delta};\mu)}{\lVert v_{\delta} \rVert_{\mathbf{V}}},\quad \forall v_{\delta} \in \mathbf{V}_{\delta}.
\end{equation}
Let us assume to have a lower bound $\alpha_{\text{LB}}(\mu)$ for the discrete coercivity constant of $a(\cdot,\cdot;\mu)$, independent on $N_{\delta}$. Let us define the following error estimator for the energy norm:
\begin{equation}
    \eta_{\text{en}}(\mu)=\frac{\lVert\hat{r}_{\delta}(\mu)\rVert_{\mathbf{V}}}{\alpha_{\text{LB}}^{1/2}(\mu)}.   
\end{equation}
Then the error between FOM and ROM solution is bounded as follows:
\begin{equation}
    \lVert u_{\delta}(\mu)-u_{\textup{rb}}(\mu) \rVert_{\mu}\le \eta_{\textup{en}}(\mu), \quad  \text{$\forall\mu \in \mathbf{P}$}.
    \label{prima}
\end{equation}
Another quantity of interest for the error analysis is the effectivity index:
\begin{equation}
    \text{eff}_{\text{en}}(\mu) = \frac{\eta_{\text{en}}(\mu)}{\lVert u_{\delta}(\mu)-u_{\text{rb}}(\mu) \rVert_{\mu}}.
    \label{eff_def}
\end{equation}
Such a quantity is $\geq 1$ and the quality of the error estimator $\eta_{\text{en}}(\mu)$ increases as the effectivity index is closer to one. If the problem is coercive and compliant, the effectivity index can be bounded as follows:
\begin{equation}
    \textup{eff}_{\textup{en}}(\mu)\le \sqrt{\gamma_{\delta}(\mu)/\alpha_{\textup{LB}}(\mu)}, \quad \text{$\forall\mu \in \mathbf{P}$},
    \label{eff_en}
\end{equation}
where $\gamma_{\delta}(\mu)$ is the discrete continuity constant of $a(\cdot,\cdot;\mu)$.
More details as well as the proofs of the inequalities \eqref{prima} and \eqref{eff_en} can be found in \cite{hesthaven2016certified,siena2023introduction}.


Eq. \eqref{varitional-problem-rb} can be formulated in matrix form: evaluates $s_{\text{rb}}(\mu) = f(u_{\text{rb}}(\mu);\mu)$, where $u_{\text{rb}}(\mu)$ satisfies
  \begin{equation}
      a(u_{\text{rb}}(\mu),\xi_j;\mu) = f(\xi_j;\mu) \qquad  1 \le j \le N.
  \end{equation}
By using eq. \eqref{varitional-problem-rb2} 
it becomes: evaluates $s_{\text{rb}}(\mu) = {(u^{\mu}_{\text{rb}})}^T f^{\mu}_{\text{rb}}$ such that
\begin{equation}
A^{\mu}_{\text{rb}}u^{\mu}_{\text{rb}}=f^{\mu}_{\text{rb}},
\end{equation}
where
\begin{equation}
    (A^{\mu}_{\text{rb}})_{j,i} =  a(\xi_i,\xi_j;\mu) \qquad \text{and} \qquad (f^{\mu}_{\text{rb}})_j = f(\xi_j;\mu), \qquad 1 \le i,j \le N.
\end{equation}
To conclude, the main benefits of projection-based ROMs regard a direct connection with the physics of the problem because they are derived from the original governing equations. This approach ensures that the ROM retains key physical features. They provide a consistent numerical solution throughout the entire simulation process and  can handle scenarios where limited data are available \cite{amsallem2012stabilization,balajewicz2012stabilization,karatzas2020projection}. 


\subsection{Data-driven ROMs}
\label{sec:data-rom}
Unlike projection-based techniques, data-driven ROMs prevent the use of the governing equations and of any physical information about the system at hand. 
Typically, the speed-up obtained with these methods is of several orders of magnitude higher than intrusive ROMs \cite{aromabook}. Thus they 
are advantageous in real-time applications where quick decision-making is crucial such as control system and optimization. 
On the other hand, the performance in terms of accuracy strongly depends on the amount and
quality of the data, and in particular on their embedded capacity to capture the underlying dynamics. Moreover, non-intrusive ROMs suffer of a complete lack of error estimation theory. 

In the context of data-driven ROMs, regression models, neural networks, or Gaussian processes, are commonly used to learn patterns and relationships between data \cite{milano2002neural,chen2021physics,hesthaven2018non,dar2023artificial}. These techniques can identify important features and correlations that might be challenging to capture with traditional modeling approaches, especially when the system exhibits a strongly non-linear behaviour \cite{papapicco2021neural,fresca2022pod,boukraichi2023parametrized}.


Even if the integration of data, machine learning methods and dimensionality reduction enables the creation of accurate and adaptable models, many challenges are still open. They include the need for representative data to train accurate models and the need to prevent overfitting \cite{daniel2020model}, especially when dealing with a limited amount of data.

Let us introduce more specifically data-driven ROMs for our purpose. 
Unlike intrusive techniques, the coefficients $u_{\text{rb}}^{\mu}$ are not computed by solving a differential algebraic system derived from the original governing equations of the problem but by computing a map in the parameter space.
Specifically, given the solution manifold (or the database of interest which can ideally include also experimental data), 
the corresponding reduced coefficients $u_{\text{rb}}^{\mu}$ are computed by projecting the high-fidelity snapshots onto the reduced space. 
The set $\{ (\mu, u_{\text{rb}}^{\mu}) : \mu \in \mathbf{P}_h \}$ is interpolated through various techniques, among Neural Networks (NNs) and Radial Basis Functions (RBFs). Now we provide some details in more about these two techniques. 
\subsubsection{Neural Networks} 
A NN is a computational model that can learn from observations. The main benefit of NNs is their \emph{universal approximators} feature \cite{cybenko1989approximation}. They are able to learn non-linear relationships via a trial and error process. This is useful especially if the ROM involves many physical and geometrical configurations (as often occurs in real-world applications). However, the main limitations regard the setting of the optimal configuration of the network and therefore the time required for repeated training as well as the amount of data needed.

A standard NN consists of a certain number of neurons and a set of directional weighted synaptic connections among them. Neurons act as nodes and synapses as edges in an oriented graph. 
Three functions identify a neuron $j$:
\begin{itemize}
    \item the propagation function $u_j$: 
    $$u_j=\sum_{k=1}^m w_{s_k,j}y_{s_k} + b_j,$$ where $b_j$ is the bias, $y_{s_k}$ is the output related to the sending neuron $k$,  $w_{s_k,j}$ are the weights and $m$ is the number of sending neurons connected with the neuron $j$. 

    \item the activation function $a_j$: 
    $$ a_j=f_{\text{act}}\left(\sum_{k=1}^m w_{s_k,j}y_{s_k}+b_j\right).$$
    Possible choices are sigmoid function, hyperbolic tangent, RELU, SoftMax \cite{sharma2017activation}. 
    \item the output function $y_j$:
    \begin{equation*}
        y_j=f_{\text{out}}(a_j). 
    \end{equation*}
    Often the output function coincides with the identity, so that $a_j = y_j$.
\end{itemize}

Typically feedforward NNs, where neurons are arranged into layers, are used for interpolation purpose \cite{rosenblatt1958perceptron,fine2006feedforward}. They consist of input and output layers, with one or more hidden layers in between. 
The weights of the network are adjusted during the training process to minimize the difference between the predicted outputs and the actual values in the training dataset. This quantity is measured by means of the loss function $\mathcal{L}$: the mean squared error (MSE) is  a common choice for $\mathcal{L}$ for interpolation problems. The backpropagation algorithm is employed \cite{rojas1996backpropagation,rumelhart1986learning} which uses the gradient of $\mathcal{L}$ with respect to the weights:
\begin{equation}
\begin{split}
    & \frac{\partial \mathcal{L}}{ \partial w_{s_k,j}^l}=\frac{\partial \mathcal{L}}{\partial a_j^l}\frac{\partial a_j^l}{\partial u_j^l}\frac{\partial u_j^l}{\partial w_{s_k,j}^l} , \\
    & \frac{\partial \mathcal{L}}{ \partial b_j^l}=
   \frac{\partial \mathcal{L}}{\partial a_j^l}\frac{\partial a_j^l}{\partial u_j^l} \frac{\partial u_j^l}{\partial b_j^l}.
\end{split}
\end{equation}
This algorithm allows to compute the gradient of the error with respect to the weights for a given input by propagating it backwards through the network.  In particular, the forward pass computes the values of the output layer from the input one and evaluates the loss function. Backpropagation performs a backward pass after each forward pass to  compute the gradient of the loss function and adjust the model parameters as follows:
\begin{equation}
\begin{split}
    &  \bm{w}=\bm{w}-\eta \frac{\partial \mathcal{L}}{\partial \bm{w}}, \\
    & \bm{b} = \bm{b} -\eta  \frac{\partial \mathcal{L}}{\partial \bm{b}},
\end{split}
\end{equation}
where $\bm{w}$ and $\bm{b}$ are matrix representations of the weights and biases and $\eta$ is the learning rate. 
Notice that tuning the hyperparameters, such as the learning rate, the activation function, and the number of hidden layers and neurons, is essential for achieving an optimal performance of the network \cite{goodfellow2016deep,kriesel2007brief,calin2020deep}. 

 \subsubsection{Radial basis function} 
 
 Interpolation by RBFs is a standard approach to approximate functions. RBFs are symmetric and depend only on the distance from a center point. Commonly used RBFs include Gaussian functions, thin-plate splines, and multiquadrics \cite{franke1982scattered,buhmann1993spectral,hardy1971multiquadric}. The main advantage of RBF interpolation method is that it provides a smooth and continuous interpolating function without the oscillations often present in polynomial-based interpolation techniques. However, finding optimal weights can be challenging, the computational cost could be larger compared to polynomial interpolation methods, and the interpolation in more variables could be affected by several issues (\cite{buhmann2000radial,forti2014efficient,skala2016practical}). 
 
 If $f(x)$ is the function to approximate with $x$ the independent variable, the RBF interpolation provides:
\begin{equation}
    f(x) = \sum_{i=1}^n w_i\phi(\| x - x_i\|) + P(x),
\end{equation}
where $n$ is the number of the interpolation points, $x_i$ are the interpolation points, $w_i$ are the weights, $\phi$ is the chosen RBF and $P(x)$ is a polynomial required for stability reasons \cite{lindner2017radial}. 
Given the input-output pairs, ${(x_i,y_i)}_{i=1}^n$, the weights are determined by solving a linear system of equations derived from the interpolation conditions:
\begin{equation}
    f(x_i) = y_i, \quad i =1,\dots,n.
\end{equation}
By adding the constrains:
\begin{align*}
	\sum_{i=1}^{n}w_{i} = 0\qquad\text{and}\qquad \sum_{i=1}^{n}w_{i}{x}_{i} = 0
\end{align*}
the weights $w_{i}$ and the polynomial $P$ are uniquely determined.


\subsection{Reduced basis generation}\label{sec:pod}

In this section we analyze two strategies for generating a reduced basis space: Proper Orthogonal Decomposition (POD) and greedy algorithm. Many other (linear and non-linear) techniques, such as Proper Generalized Decomposition (PGD), factor analysis, independent component analysis and autoencoders can be found in the literature \cite{cunningham2015linear,fu2021data,lee2020model,gonzalez2018deep,kashima2016nonlinear}. 

Starting from an equispaced or a random sampling of $\mathbf{P}$, let us introduce a finite-dimensional set of parameters $\mathbf{P}_{h}$. 
Then the following high-fidelity solutions manifold can be provided:
\begin{equation}
    \mathcal{M}_{\delta}(\mathbf{P}_h)=\{ u_{\delta}(\mu) : \mu \in \mathbf{P}_h  \},
\end{equation}
whose dimension is $M=|\mathbf{P}_h|$. 
By definition we have:
\begin{equation}
    \mathcal{M}_{\delta}(\mathbf{P}_h)\subset\mathcal{M}_{\delta} \qquad \text{and} \qquad \mathbf{P}_h \subset \mathbf{P}.
\end{equation}
Therefore, if $\mathbf{P}_h$ is rich enough, $\mathcal{M}_{\delta}(\mathbf{P}_h)$ results to be a proper approximation of $\mathcal{M}_{\delta}$. 

\subsubsection{Proper orthogonal decomposition}

POD is a method suitable for data compression \cite{eckart1936approximation}. Once the parameter space is discretized and FOM solutions are collected for each element of $\mathbf{P}_h$, POD allows to retain only the essential information about the system at hand. 
The $N$-dimensional POD space can be computed as the solution to the following minimization problem:
\begin{equation}
    \min_{\mathbf{V}_{\text{rb}}:|\mathbf{V}_{\text{rb}}|=N} \Bigg( \int_{\mu \in \mathbf{P}} \inf_{v_{\text{rb}}\in \mathbf{V}_{\text{rb}}} \lVert u_{\delta}(\mu)-v_{\text{rb}} \rVert_{\mathbf{V}}^2 d\mu \Bigg)^{1/2}. \label{eq:POD1}
\end{equation}
In discrete form eq. \eqref{eq:POD1} yields:
\begin{equation}
     \min_{\mathbf{V}_{\text{rb}}:|\mathbf{V}_{\text{rb}}|=N} \Bigg( \frac{1}{M}\sum_{\mu \in \mathbf{P}_h} \inf_{v_{\text{rb}}\in \mathbf{V}_{\text{rb}}} \lVert u_{\delta}(\mu)-v_{\text{rb}} \rVert_{\mathbf{V}}^2  \Bigg)^{1/2}.
     \label{opt-crit}
\end{equation}
Let the elements of $\mathbf{P}_h$ be 
\begin{equation}
    \{\mu_1,\dots,\mu_M\},
\end{equation}
and indicate with
\begin{equation}
    \{ \psi_1, \dots, \psi_M \}
\end{equation}
the elements of $\mathcal{M}_{\delta}(\mathbf{P}_h)$, $\psi_m=u_{\delta}(\mu_m)$ for $m=1,\dots,M$ (i.e. the high-fidelity snapshots). Furthermore, let define $\mathbf{V}_{\mathcal{M}}=\text{span}\{ u_{\delta}(\mu) : \mu \in \mathbf{P}_h \}$ and the symmetric and linear operator $\mathcal{C}:\mathbf{V}_{\mathcal{M}}\mapsto\mathbf{V}_{\mathcal{M}}$ such that:
\begin{equation}
    \mathcal{C}(v_{\delta})=\frac{1}{M}\sum_{m=1}^M (v_{\delta},\psi_m)_{\mathbf{V}} \psi_m, \quad v_{\delta} \in \mathbf{V}_{\mathcal{M}}.
    \label{linear-op}
\end{equation}
Its eigenvalues and normalized eigenvectors are identified as $(\lambda_i,\xi_i)\in \mathbf{R}\times \mathbf{\mathbf{V}_{\mathcal{M}}}$ and fulfill:
\begin{equation}
    (\mathcal{C}(\xi_i),\psi_m)_{\mathbf{V}}=\lambda_i(\xi_i,\psi_m)_{\mathbf{V}}, \quad 1 \le m \le M,
    \label{eig-p}
\end{equation}
with 
$$\lambda_1 \ge \lambda_2 \ge \dots \ge \lambda_M.$$
The eigenvectors $\{\xi_1,\dots,\xi_M\}$ represent the basis functions for the space $\mathbf{V}_{\mathcal{M}}$ and the first $N \ll M$ eigenfunctions $\{\xi_1,\dots,\xi_N\}$ give the $N$-dimensional reduced space
\begin{equation}
    \mathbf{V}_{\text{POD}} = \text{span} \{\xi_1,\dots,\xi_N \}.
\end{equation}

The error introduced by replacing the elements of $\mathcal{M}_{\delta}(\mathbf{P}_h)$ with $\mathbf{V}_{\text{POD}}$  is given by the sum of the neglected eigenvalues \cite{quarteroni2015reduced}:
\begin{equation}
    \frac{1}{M}\sum_{m=1}^M\lVert \psi_m - P_N[\psi_m] \rVert_{\mathbf{V}}^2 = \sum_{m=N+1}^M \lambda_m,
\end{equation}
where $P_N[\psi_m] = \sum_{i=1}^N(\psi_m,\xi_i)_{\mathbf{V}}\xi_i$ is the projection of $\psi_m$ onto $\mathbf{V}_{\text{POD}}$.
As a consequence of the orthonormality of the eigenvectors in $\lVert \cdot \rVert_{\ell^2(\mathbf{R}^M)}$, we have:
\begin{equation}
    (\xi_m,\xi_q)_{\mathbf{V}}=M\lambda_i\delta_{mq}, \quad 1\le m,q \le M,
\end{equation}
where $\delta_{mq}$ is the Kronecker delta.

The correlation matrix $C\in\mathbf{R}^{M\times M}$ representing the linear operator \eqref{linear-op} can be expressed as:
\begin{equation}
    C_{m,q}=\frac{1}{M}(\psi_m,\psi_q)_{\mathbf{V}}, \quad 1 \le m,q \le M.
\end{equation}
Therefore the eigenvalues problem \eqref{eig-p} is reformulated as:
\begin{equation}
    C v_i = \lambda_i v_i, \quad 1 \le i \le N.
\end{equation}
Finally, the orthogonal basis functions are given by:
\begin{equation}
    \xi_i = \frac{1}{\sqrt{M}}\sum_{m=1}^M (v_i)_m \psi_m, \quad 1 \le i \le N,
\end{equation}
where$(v_i)_m$ denotes the $m$-th element of the eigenvector $v_i \in \mathbf{R}^M$ .

The computational cost to implement the POD algorithm can be very large, because it is impossible to know how many high-fidelity solutions are needed to guarantee a good approximation of the system behaviour and it depends on the problem at hand. Typically $M\gg N$ instances of the FOM need to be solved during the offline stage and, when $M$ and $N_{\delta}$ are large, the cost of extracting the eigenfunctions rises, as it scales as $\mathcal{O}(NN_{\delta}^2)$.

\subsubsection{Greedy algorithm} 
The greedy algorithm is an iterative method which add at each iteration a basis function, until the desired accuracy is achieved \cite{bang2004greedy}. A detailed scheme is shown in Algorithm \ref{alg:the_alg}.
This strategy needs to compute one FOM solution per iteration and therefore $N$ FOM solutions to build the $N$-dimensional reduced basis space. 

Let us introduce an estimation $\eta(\mu)$ of the error between the FOM the ROM solution (the computation of this bound is out of the scope of the chapter, more details can be found in \cite{hesthaven2016certified,siena2023introduction}):
\begin{equation}
    \lVert u_{\delta}(\mu)-u_{\text{rb}}(\mu) \rVert_{\mu} \le \eta(\mu), \quad \forall \mu \in \mathbf{P}.
\end{equation}
At the $n$-th step of the iterative procedure, the parameter $\mu_{n+1}$ is chosen:
\begin{equation}
    \mu_{n+1}=\underset{\mu \in \mathbf{P}}{\operatorname{argmax}} \big[  \eta(\mu) \big],
\end{equation}
and the related high fidelity solution $u_{\delta}(\mu_{n+1})$ is stored in the reduced basis space $\mathbf{V}_{\text{rb}}=\text{span} \{u_{\delta}(\mu_1),\allowbreak\dots,\allowbreak,u_{\delta}(\mu_{n+1}) \} $.
It means that at each iteration the $n+1$ basis function optimizes the approximation of the ROM error over the parameter set $\mathbf{P}$, until the maximum error falls below a specified tolerance. Therefore, while the greedy method relies on the maximum norm over $\mathbf{P}$, the POD algorithm looks for an optimal basis in the $L^2$-norm. 

From the operational standpoint, the finite dimensional parameter set $\mathbf{P}_h$ is provided as for the POD approach. However, the greedy method does not require the resolution of the problem \eqref{varitional-problem-discr} for each point in $\mathbf{P}_h$, but it only needs an estimation for $\eta(\mu)$. Therefore, the parameter space can be richer
than the one used with the POD method, just because the computational cost is significantly lower. So the main benefits of this approach, with respect to POD, are to prevent the resolution of a big eigenvalues problem and to compute a significantly lower number of FOM
solutions.

More in general, 
when $F=\{ f(\mu):\Omega \mapsto \mathbf{R} : \mu \in \mathbf{P} \}$ is a collection of parametrized functions, the $(n+1)$-th basis function is computed as
\begin{equation}
    f_{n+1} = \underset{\mu \in \mathbf{P}}{\operatorname{argmax}}  \lVert f(\mu)- P_n f(\mu) \rVert_{\mathbf{V}} ,
\end{equation}
where $P_nf$ is the orthogonal projection of $f$ onto $F_n=\text{span}\{f_1,\dots,f_n \}$. 
The convergence of the greedy method is ensured when $F$  has an exponentially small Kolmogorov $N$-width. In that case, the reduced basis approximation $u_{\text{rb}}(\mu)$ goes exponentially to the high fidelity solution $u_{\delta}(\mu)$ (more details about the convergence are analyzed in \cite{hesthaven2016certified,siena2023introduction}). Remember that it is not enough to declare that the ROM works properly, because it is essential to assume that FOM solution is a proper estimation of the exact solution $u(\mu)$. 
Finally, it is crucial to remark that the snapshots computed during the greedy procedure for each parameter in $\mathbf{P}_h$, 
may have a linear relationship, resulting in a strongly ill conditioned solution matrix. Therefore, in order to solve this problem, the snapshots $ u_{\delta}(\mu_1),\dots, u_{\delta}(\mu_N) $ should be orthonormalized to compute the basis $\xi_1,\dots,\xi_N$.

\begin{algorithm}
\caption{The greedy algorithm}\label{pod-alg}
\text{\textbf{Input:} $\mu_1$, tol}\\
\text{\textbf{Output:} $\mathbf{V}_{\text{rb}}$}
\begin{algorithmic}[1]
\State Compute $u_{\delta}(\mu_n)$
\State $\mathbf{V}_{\text{rb}}=\text{span}\{u_{\delta}(\mu_1),\dots,u_{\delta}(\mu_n) \}$
 \For{$\mu \in \mathbf{P}_h$}
 \State compute $u_{\text{rb}}(\mu)$
  \State evaluate $\eta(\mu)$
 \EndFor
\State $ \mu_{n+1}=\underset{\mu \in \mathbf{P}}{\operatorname{argmax}} \big[  \eta(\mu) \big],$

\If{$\eta(\mu_{n+1}) >$ tol}
    \State $n=n+1$

\Else{}
    \State break.
\EndIf

\end{algorithmic}
\label{alg:the_alg}
\end{algorithm}

\section{Test cases}\label{sec:app}
In this section we  test the performance in terms of efficiency and accuracy of ROM approaches described in Sec. \ref{sec:rom}. 
We start with an academic benchmark, the thermal block problem (Sec. \ref{sec:termal}). In such a case, we adopt a Galerkin projection ROM method where the basis functions are computed by using the greedy algorithm. Then we move to two real-world applications, the one related to the analysis of the blood flow patterns in coronary arteries (Sec. \ref{sec:coronary}) and the other one related to the modelling of the granulation process in pharmaceutics industry (Sec. \ref{sec:industry}). Here the reduced space is extracted by POD and the reduced coefficients are computed by using a data-driven approach, based on NNs for the biomedical problem and based on RBFs for the industrial one.

\subsection{Thermal block problem}\label{sec:termal} 

\begin{figure}
\centering
\begin{tikzpicture}[thick,scale=5]

\draw[black] (0.5,0.5) -- (-0.5,0.5);

\draw[teal] (0.5,-0.5) -- (0.5,0.5);
\draw[teal] (-0.5,0.5) -- (-0.5,-0.5);

\draw[brown] (-0.5,-0.5) -- (0.5,-0.5);

\draw[fill=red!10] (0,0) circle [radius=0.25cm];
\node[color=black] at (0,0) {$\Omega_1$};
\node at (-0.4,0.4) {$\Omega_2$};

\node at (0.6,0) {\textcolor{teal}{$\Gamma_{\text{side}}$}};

\node at (0,-0.55) {\textcolor{brown}{$\Gamma_{\text{base}}$}};

\node at (0,0.55) {\textcolor{black}{$\Gamma_{\text{top}}$}};

\end{tikzpicture}
\caption{Thermal block problem: computational domain.}
\label{fig:tblock}
\end{figure}
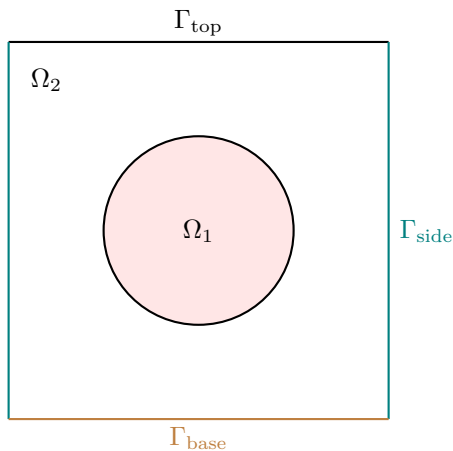

Let us consider a bi-dimensional steady heat conduction problem over the domain $\Omega=(-1,1)\times(-1,1)$ shown in Figure \ref{fig:tblock}. 
The subdomain $\Omega_1$ is a disk centered at the origin of the domain with radius $r_0=0.5$ such that  $\Omega = \Omega_1 \cup \Omega_2$ and $\Omega_1 \cap \Omega_2 = \emptyset$. The boundary $\partial\Omega$ is composed by $\Gamma_{\text{side}}\cup\Gamma_{\text{top}}\cup\Gamma_{\text{base}}$ as shown in Figure \ref{fig:tblock}. In particular, $\Gamma_{\text{base}}=(-1,1)\times\{-1\}$, $\Gamma_{\text{top}}=(-1,1)\times\{1\}$ and $\Gamma_{\text{side}}=\{\pm1\}\times(-1,1)$.

The thermal conductivity $k$ is constant on each subdomain:
\begin{equation}
    k = \begin{cases}
        k_0 & \mbox{on } \Omega_1, \\
        1 & \mbox{on } \Omega_2. \\
    \end{cases}
\end{equation}
The parameters considered for this problem are $\mu_0=k_0$ and the heat flux $\mu_1$ through $\Gamma_{\text{base}}$. The parameter vector $ \mu = (\mu_0, \mu_1) \in \mathbf{P}$ ranges in $\left[0.1, 10\right]\times\left[-1, 1\right]$.  
Concerning other boundaries, $\Gamma_{\text{side}}$ is thermally isolated and on $\Gamma_{\text{top}}$ a reference temperature is employed. 

Therefore, the FOM in strong formulation is given by 
\begin{equation}
    \begin{cases}
        -\nabla \cdot (k\nabla u(\mu)) = 0, & \mbox{on }\Omega,\\
        u(\mu) = 0, & \mbox{on }\Gamma_{\text{top}},\\
        k\nabla u(\mu) \cdot \bm n = 0, & \mbox{on }\Gamma_{\text{side}},\\
        k\nabla u(\mu) \cdot \bm n = \mu_1, & \mbox{on }\Gamma_{\text{base}},
    \end{cases}
    \label{eq:fomtblock}
\end{equation}
where $u(\mu)$ is the temperature and $\bm n$ is the outward normal to the boundary.

Let us consider the Hilbert space $\mathbf{V}=\{v\in H^1(\Omega):v_{\mid \Gamma_{\text{top}}}=0\}$. Then the weak formulation of the problem \eqref{eq:fomtblock} evaluates $s(\mu) = l(u;\mu)$ where $u(\mu)\in \mathbf{V}$ satisfies:
\begin{equation}
      a(u,v;\mu) = f(v;\mu) \qquad \forall v \in \mathbf{V}.
\end{equation}
The bilinear form is $a(u,v;\mu) = \int_{\Omega} k \nabla u \cdot \nabla v d\Omega$ while the linear forms are $f(v; \mu)=\mu_1\int_{\Gamma_{\text{base}}} v dA$ and $l(u;\mu)=\int_{\Gamma_{\text{base}}} u(\mu)dA$ with $dA$ an infinitesimal portion of the boundary. In our case, the output of interest $s(\mu)$ is the average temperature over $\Gamma_{\text{base}}$. Note that the problem is compliant and affine, in fact 
\textcolor{black}{it can be written as$$ \mu_0{\int_{\Omega_1} \nabla u \cdot \nabla v d\Omega} + {1} \int_{\Omega_2} \nabla u \cdot \nabla v d\Omega ={\mu_1}{\int_{\Gamma_{\text{base}}} v dA},$$ therefore the parameter-dependent terms in \eqref{affine} are $$\theta^0_a = \mu_0, \quad \theta^1_a = 1,  \quad \theta^0_f = \mu_1,$$ and the parameter-independent terms in \eqref{affine} are $$a_0 = \int_{\Omega_1} \nabla u \cdot \nabla v d\Omega, \quad a_1 = \int_{\Omega_2} \nabla u \cdot \nabla v d\Omega, \quad f_0 = \int_{\Gamma_{\text{base}}} v dA.$$}


For the well-posedness of the problem, we assume that $\mu_0>0$, therefore $k>0$ and the coercivity of $a(\cdot,\cdot;\mu)$ is verified. The continuity
of $a(\cdot,\cdot;\mu)$ and $f(\cdot;\mu)$ can be obtained with the Cauchy-Schwarz inequality. 
Therefore, the Lax-Milgram theorem guarantees the existence and uniqueness of the solution. A finite element method employing piece-wise linear elements is adopted to compute FOM solutions \cite{hesthaven2016certified}. The mesh consists of 542 cells. 

Here we use a certified reduced basis approximation based on the greedy algorithm (see Sec. \ref{sec:pod}). 
We  consider 100 discrete values 
of the parameters to train the ROM and the same number for the validation. The reduced basis are computed with a tolerance of $10^{-5}$ corresponding to a reduced dimension of $N=4$. An error analysis is performed varying the number of basis functions until $N=4$. Hereinafter, unless otherwise specified, the error we are referring is 
\begin{equation*}
\|u_{\delta}(\mu) - u_{\text{rb}}(\mu)\|_{\mathbf{V}},
\end{equation*} 
and corresponding average and normalized values.

\begin{figure}
	\centering
 	\subfloat[][Mean and maximum error for $u$.\label{fig:error_u}]{\includegraphics[height=.35\textwidth]{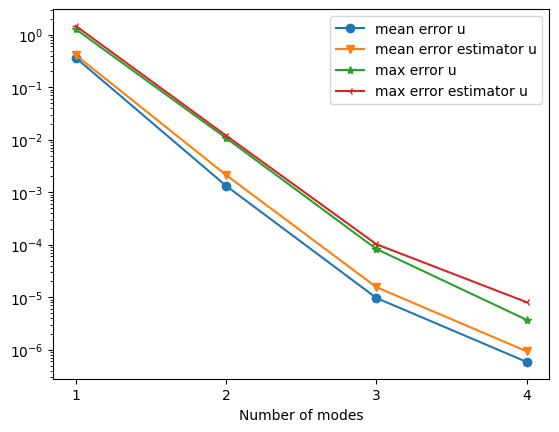}}
	\subfloat[][Mean and maximum for the effectivity index of $u$.\label{fig:effectivity_u}]{\includegraphics[height=.35\textwidth]{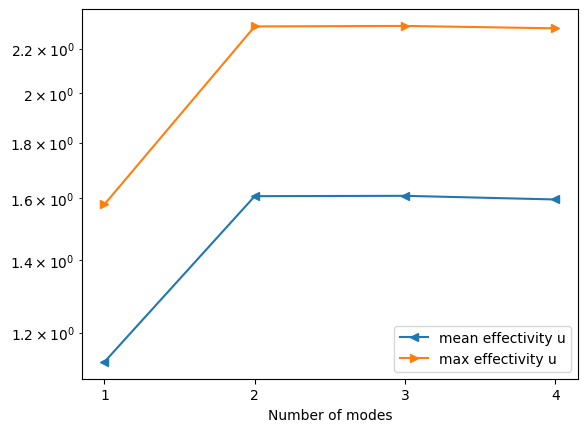}}
 
  	\subfloat[][Mean and maximum error for $s$.\label{fig:error_l}]{\includegraphics[height=.35\textwidth]{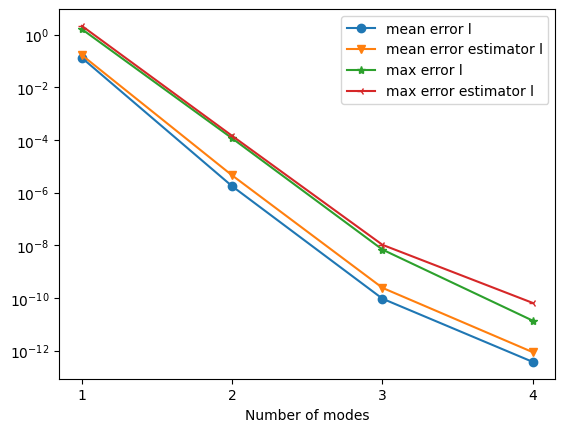}}
	\subfloat[][Mean and maximum for the effectivity index of $s$.\label{fig:effectivity_l}]{\includegraphics[height=.35\textwidth]{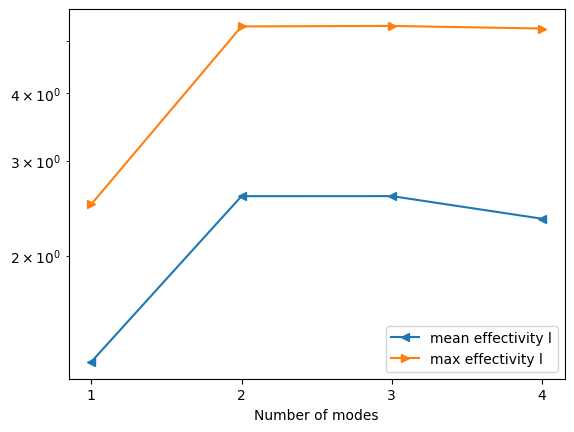}}

  	\subfloat[][Mean and maximum relative error for $u$.\label{fig:rel_error_u}]{\includegraphics[height=.35\textwidth]{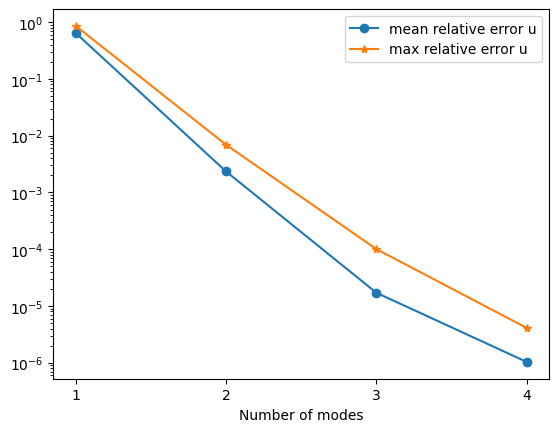}}
	\subfloat[][Mean and maximum relative error for $s$.\label{fig:rel_err_l}]{\includegraphics[height=.35\textwidth]{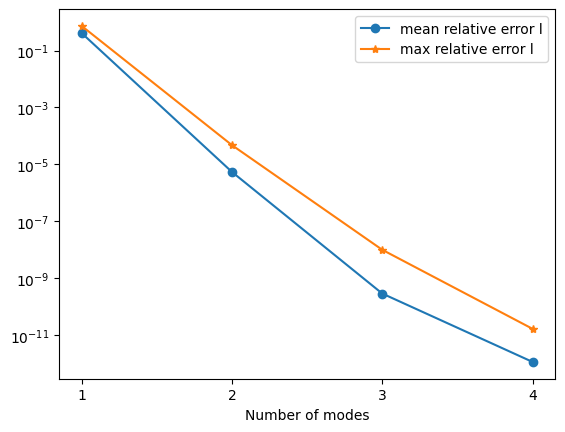}}
	\caption{Thermal block problem: error analysis for the temperature $u$ and the output of interest $s$ at varying of the number of basis functions. 
 }
	\label{eig-and-err}
\end{figure}

In Figure \ref{fig:error_u} and \ref{fig:error_l} the error is shown at varying of the number of basis functions. As expected, the error decreases as the number of modes increases. In addition, according to equation \eqref{prima}, both mean and maximum errors are bounded from the error estimators $\eta_{\text{en}}$. Consistently, the effectivity index \eqref{eff_def} in Figure \ref{fig:effectivity_u} and \ref{fig:effectivity_l} is close to 1 and it does not show large variability as the number of modes changes (both for the temperature and for the output of interest). It is relevant to observe that also the relative errors decrease with the number of basis functions (Figure \ref{fig:rel_error_u} and \ref{fig:rel_err_l}). In particular, they reach with 4 modes about the same order of magnitude of the absolute errors, respectively $10^{-6}$ for the temperature $u$ and $10^{-11}$ for the output $s$. 

\begin{figure}
	\centering
 	\subfloat[][FOM solution for $u$.\label{fig:fomt}]{\includegraphics[height=.35\textwidth]{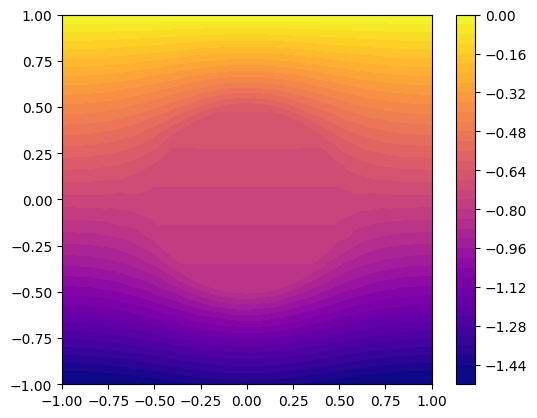}}
	\subfloat[][ROM solution for $u$. \label{fig:romt}]{\includegraphics[height=.35\textwidth]{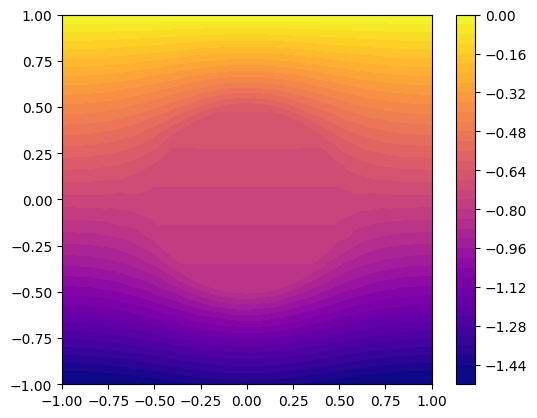}}
 \caption{Thermal block problem: qualitative comparison between FOM and ROM solution for the test point $\mu=(8,-1)$.}
 \label{fig:qualitative_thermal}
 \end{figure}

For completeness, a qualitative comparison of $u$ is show in Figure \ref{fig:qualitative_thermal} for the parameter $\bm \mu = (8,-1)$ and we can appreciate a good matching between the FOM and ROM solutions. 
The speed-up realized through this approach is around 9, which is notably positive. Nevertheless, it is crucial to recognize that these outcomes stem from a toy problem with a restricted number of degrees of freedom.

\subsection{A cardiovascular application}
\label{sec:coronary}
\textcolor{black}{Dealing with biomedical applications is complex due to uncertainties in model formulations and parameter calibration, alongside inherent variability in biological systems. These uncertainties arise from limited experimental data, noise in measurements, and the intricate nature of physiological processes. Addressing these challenges requires advanced computational techniques and interdisciplinary collaboration to ensure accurate and reliable biomedical modeling outcomes.}

In the framework of the cardiovascular applications, our aim is to use ROMs to provide accurate predictions of the blood flow in patient-specific configurations. 
In fact, recently the medical community increased demand for quantitative investigations of cardiovascular physiology and pathology, in order to act more consciously \cite{africa2024lifex_cfd, africa2023lifex_ep, africa2023epmf, africa2023lifex_fibers, fedele2023, bucelli2023, africa2022lifex, vergara2022, zingaro2022, regazzoni2022, piersanti2021, salvador2021, stella2020}. The ambition is to design a virtual platform where surgeons can have real-time outcomes about the hemodynamics indices of the patient at hand. See also Sec. \ref{sec:atlas}.

\begin{figure}
\centering

\subfloat[][Computational domain.\label{fig:cabg-domain}]{
\begin{overpic}[height=0.45\textwidth, trim = 0cm -9cm 0cm 0cm]{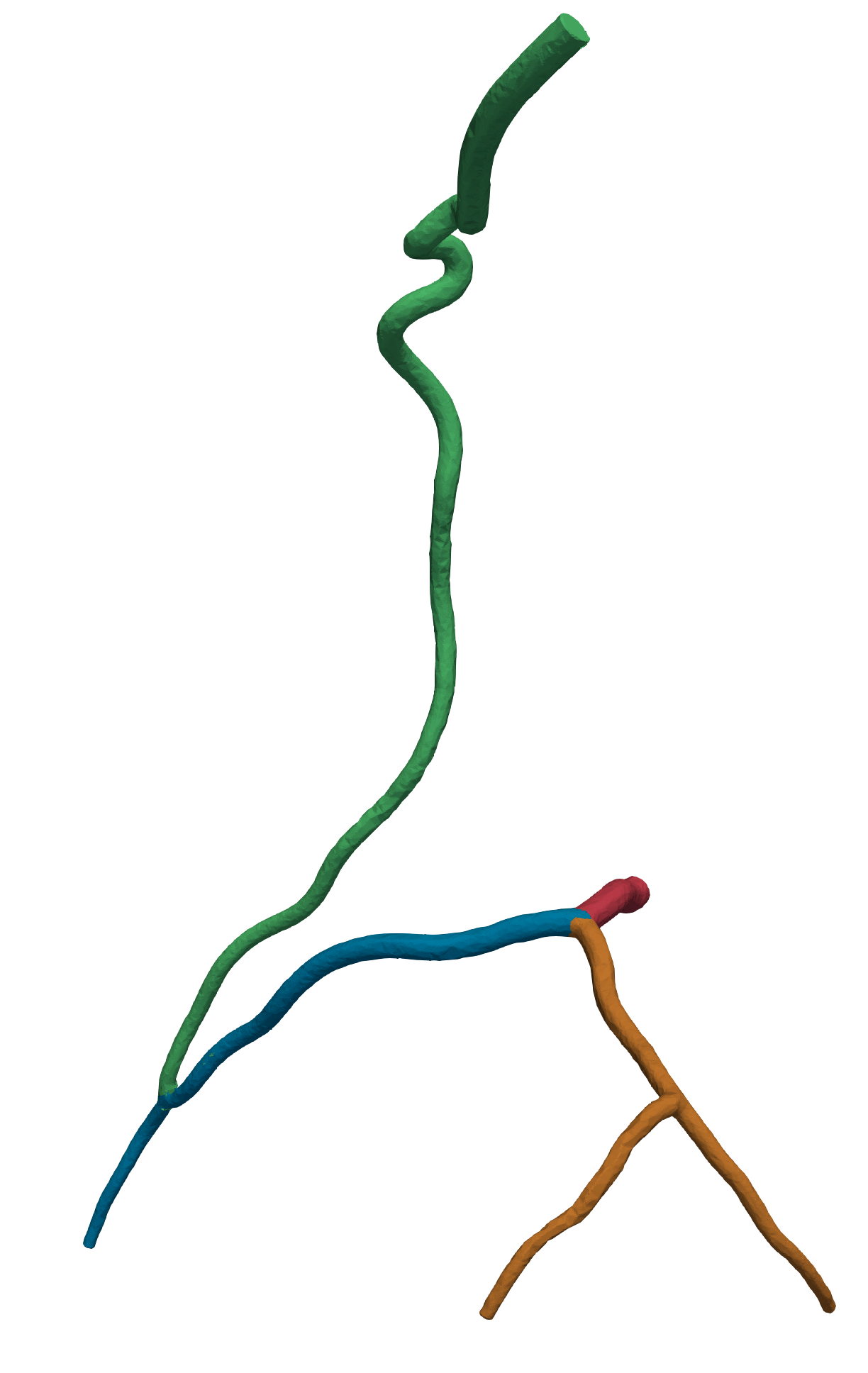}
        \put(5,80){\textcolor{emerald}{LITA}}
        \put(40,50){\textcolor{BrickRed}{LMCA}}
        \put(50,25){\textcolor{bronze}{LCx}}
        \put(12,25){\textcolor{NavyBlue}{LAD}}
        \linethickness{1pt}\put(39,107){\color{Green}\vector(-0.08,-0.22){3}}
        \linethickness{1pt}\put(48,49){\color{Green}\vector(-0.60,-0.30){7.5}}
        \linethickness{1pt}\put(5,22){\color{red}\vector(-0.08,-0.22){3}}
        \linethickness{1pt}\put(30,17.5){\color{red}\vector(-0.08,-0.22){3}}
        \linethickness{1pt}\put(52,18){\color{red}\vector(0.08,-0.22){3}}
\end{overpic}}
\hfill
\subfloat[][Boundary conditions. \label{fig:BCLITALMCA}]{
\begin{overpic}[height=0.45\textwidth]{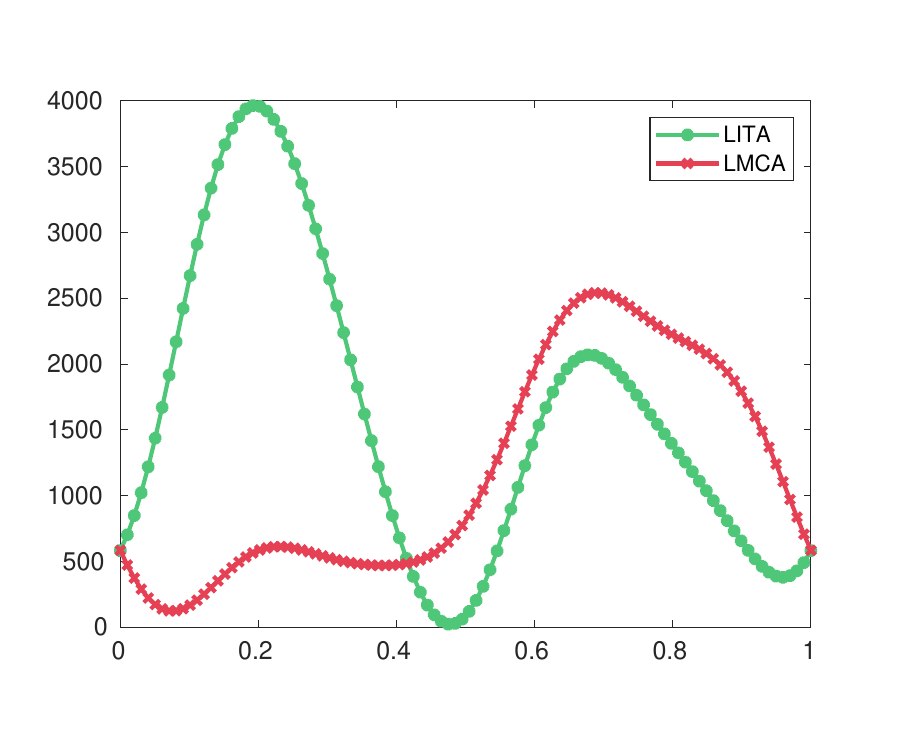}
\put(46,5){Time}
\put(-4,32){\rotatebox{90}{$\bar{q}(t)[mm^3/s]$}}
\end{overpic}}
\caption{Cardiovascular benchmark. Sketch of the computational domain: Coronary Artery Bypass Graft (CABG) with Left Internal Thoracic Artery (LITA), Left Main Coronary Artery (LMCA), Left Anterior Descending Artery (LAD) and Left Circumflex artery (LCx) (left panel). Inflow boundary conditions imposed on LITA and LMCA sections (right panel).}
\label{fig:cabg-domain-BC}
\end{figure}
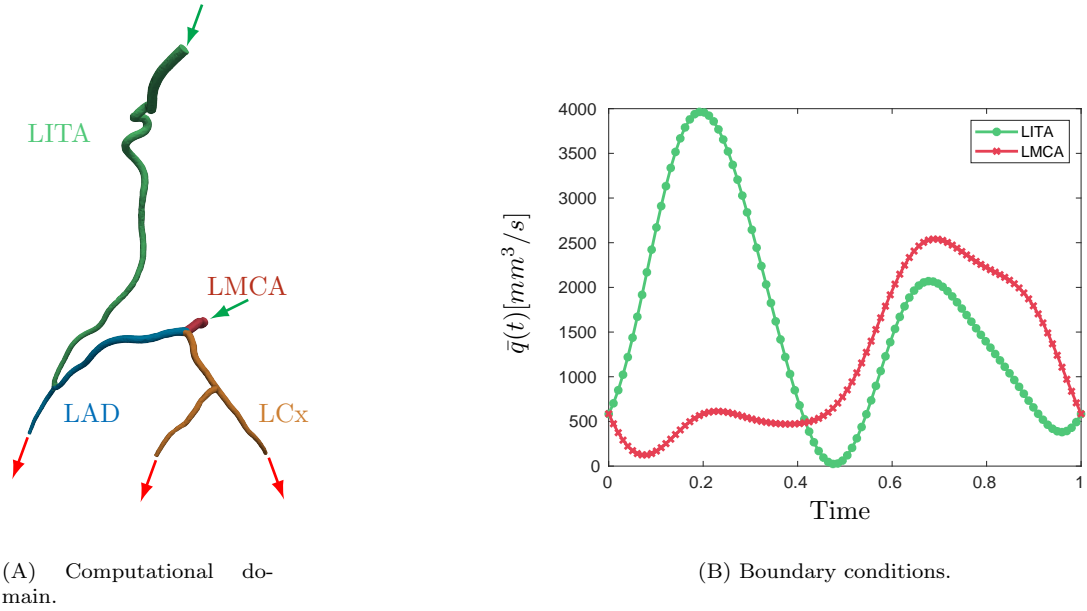

We consider a Coronary Artery Bypass Graft (CABG) when an isolated stenosis of the Left Main Coronary Artery
(LMCA) occurs. The geometry have been provided by Ospedale Luigi Sacco in Milan and it is shown in Figure \ref{fig:cabg-domain}. The model includes, beyond
the LMCA, the left internal thoracic artery (LITA), the left anterior descending artery (LAD) and the left circumflex
artery (LCx). 

The FOM in strong formulation is represented by the incompressible unsteady Navier-Stokes equations:
\begin{equation}
	\begin{cases} 
	\partial_t \bm u + \nabla \cdot (\bm u \otimes \bm u) + \nabla p - \nu \Delta \bm{u} = 0, & \mbox{on }\Omega \times (t_0, T],\\ 
	\nabla \cdot \bm u = 0, & \mbox{on }\Omega \times (t_0, T],
	\end{cases}  \label{N-S-2}
\end{equation}
where $(t_0, T]$ is the time interval of interest (with $t_0$ the initial time instant and $T$ the final time), $\bm u$ is the velocity, $p$ the kinematic pressure and $\nu$ the kinematic viscosity. A quantity of interest is the Wall Shear Stress (WSS) defined as follows:
\begin{equation}
\text{WSS} = \tau \cdot \bm{n},
\end{equation}
where $\tau = \nu (\nabla \bm{u} + \nabla \bm{u}^T)$ is the stress tensor and $\bm{n}$ is the unit normal outward vector.

The boundary $\partial \Omega$ of our computational domain includes:
\begin{itemize}
    \item two inflow, the LMCA and the LITA sections, where a realistic flow rate waveform is enforced \cite{keegan2004spiral,ishida2001mr}:
\begin{equation}
    q_i(t)=f^i\bar{q}_i(t), \quad \quad i=\text{LMCA},\text{LITA}.
    \label{flowrate}
\end{equation}
We set $f_{\text{LMCA}} = 1.12 $ and $f_{\text{LITA}}=0.82$ adapted from \cite{keegan2004spiral,ishida2001mr,verim2015cross}). The functions $\bar{q}_i(t)$ are shown in Figure \ref{fig:BCLITALMCA};
    \item the wall vessel, where no slip conditions are enforced;
    \item two outflow, the LAD and LCx sections, where  we enforce homogeneous Neumann boundary conditions. 
\end{itemize}

Unlike what was done in Sec. \ref{sec:termal}, here a second-order finite volume method is adopted for the space discretization. For sake of brevity we do not report any detail about that. The reader is referred to \cite{Siena2022,siena2023fast}. Since the problem is unsteady, we also need to introduce a time discretization. The time interval of interest $(t_0, T]$ is partitioned using the time step $\Delta t \in \mathbb{R}^+$, resulting in time levels $t^n = t_0 + n \Delta t$, where $n$ ranges from 0 to the total number of time steps $N_{T}$. Eq. \eqref{N-S-2} is discretized in time employing a Backward
Differentiation Formula of order 2 (BDF2). On the other hand, all terms in~\eqref{N-S-2} are treated implicitly except for the divergence term, which was discretized using a semi-implicit scheme. 

A data-driven ROM is employed in this case based on a POD-NN approach \cite{chen2021physics,hesthaven2018non,wang2019non} discussed in Secs. \ref{sec:data-rom} and \ref{sec:pod}. We consider a geometrical parametrization where the geometrical parameter is the degree of the stenosis in the LMCA, from mild to severe.

\begin{figure}
	\centering
 	\subfloat[][Lattice.\label{fig:lattice}]{\includegraphics[height=.35\textwidth]{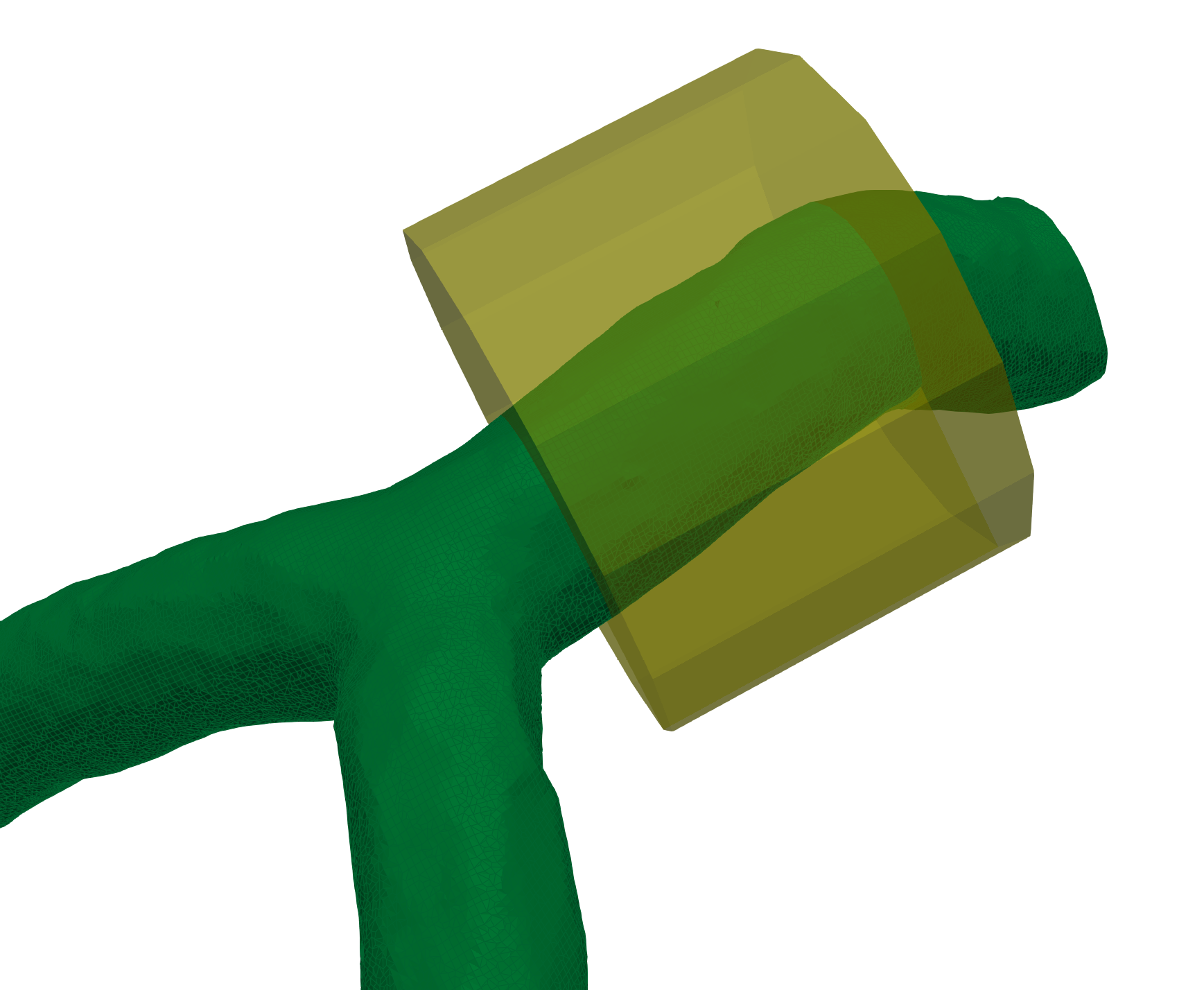}}
	\subfloat[][Deformed lattice. \label{fig:lattice_def}]{\includegraphics[height=.35\textwidth]{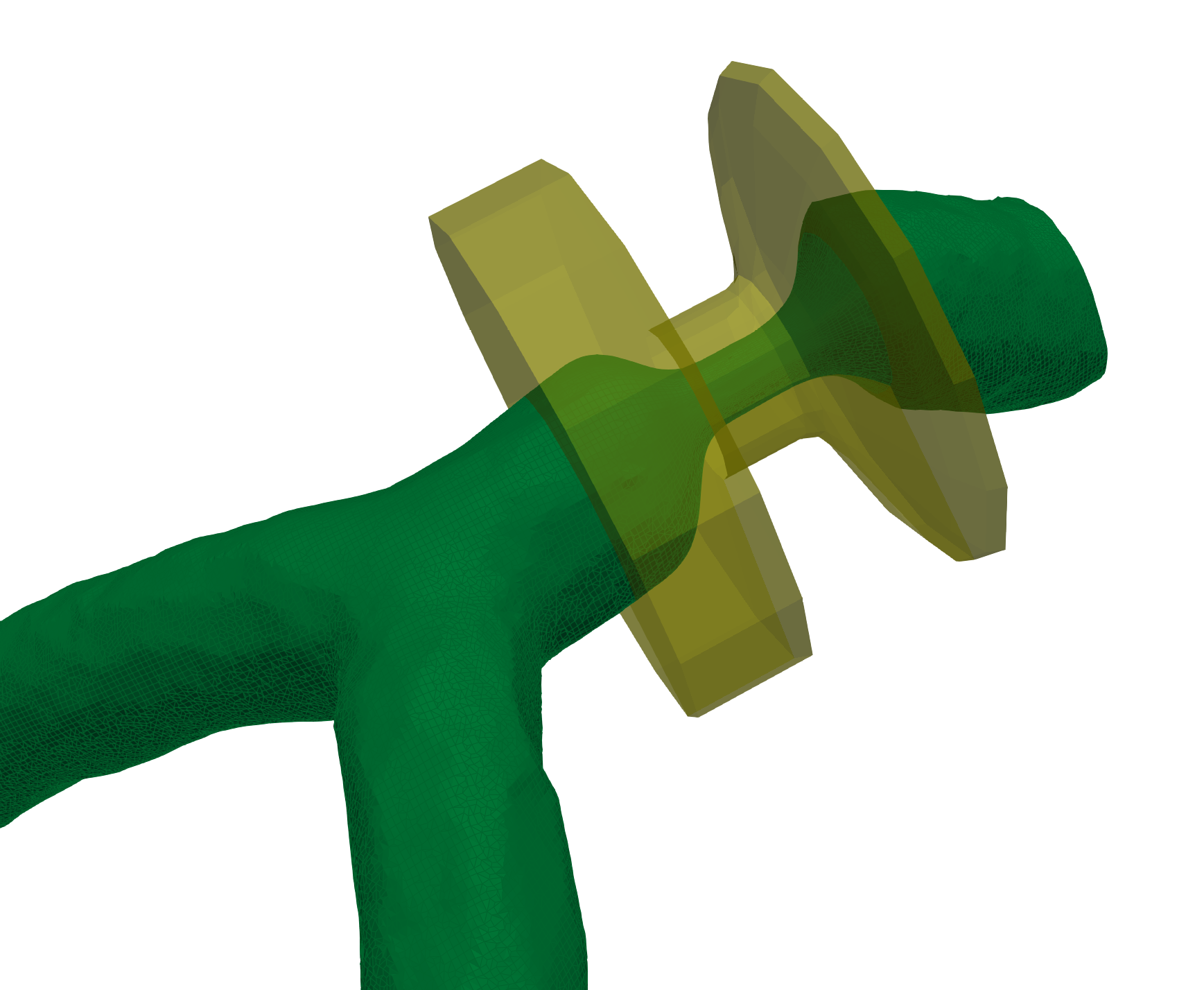}}
 \caption{Cardiovascular benchmark: introduction of stenosis in the LMCA using FFD technique.} 
 \label{fig:lattice-lattice-def}
 \end{figure}

The features of the computational mesh are shown in Table \ref{grid}.
\begin{table}
\centering
\caption{Features of the mesh. }
\label{grid}       
\begin{tabular}{cccc}
\hline\noalign{\smallskip}
Number of cells & Min/max mesh size [$m$] & Average non-orthogonality [$^\circ$] & Max skewness  \\
\noalign{\smallskip}\hline\noalign{\smallskip}
986.278 & 
3.0e-5 - 4.3e-4
& 12.9 & 2.95 \\
\noalign{\smallskip}\hline
\end{tabular}
\end{table}
To consider several stenosis severity in the LMCA compatibly with the data-driven ROM adopted, it is fundamental to warp directly the mesh and not just the geometry, so that all the deformed configurations have the same number of control volumes. At this aim the Free Form Deformation (FFD) is employed \cite{brujic2002measurement}. The control points of a parametric lattice are used to define the volume containing the portion to be warped (Figure \ref{fig:lattice}). The octree algorithm \cite{amoiralis2008freeform,lamousin1994nurbs,africa2023quadtree} is employed to find a match between the control points of the lattice and the computational domain. Then, the coordinates of the control points are modified, so that the parametric volume and consequently
the computational domain are deformed (Figure \ref{fig:lattice_def}). 
 
\begin{figure}
	\centering
    \subfloat[][Velocity. \label{fig:err-cabg-u}]{
    \begin{overpic}[height=0.35\textwidth]{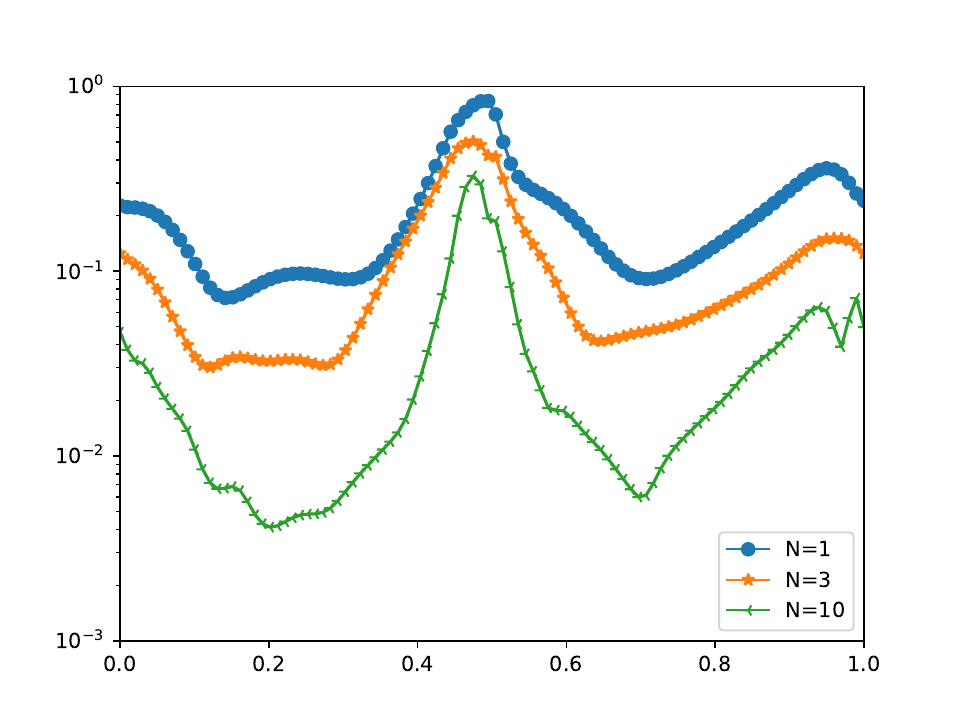}
    \put(47,0){Time}
    \put(-4,25){\rotatebox{90}{Relative error}}
    \end{overpic}}
    \subfloat[][WSS. \label{fig:err-cabg-wss}]{
    \begin{overpic}[height=0.35\textwidth]{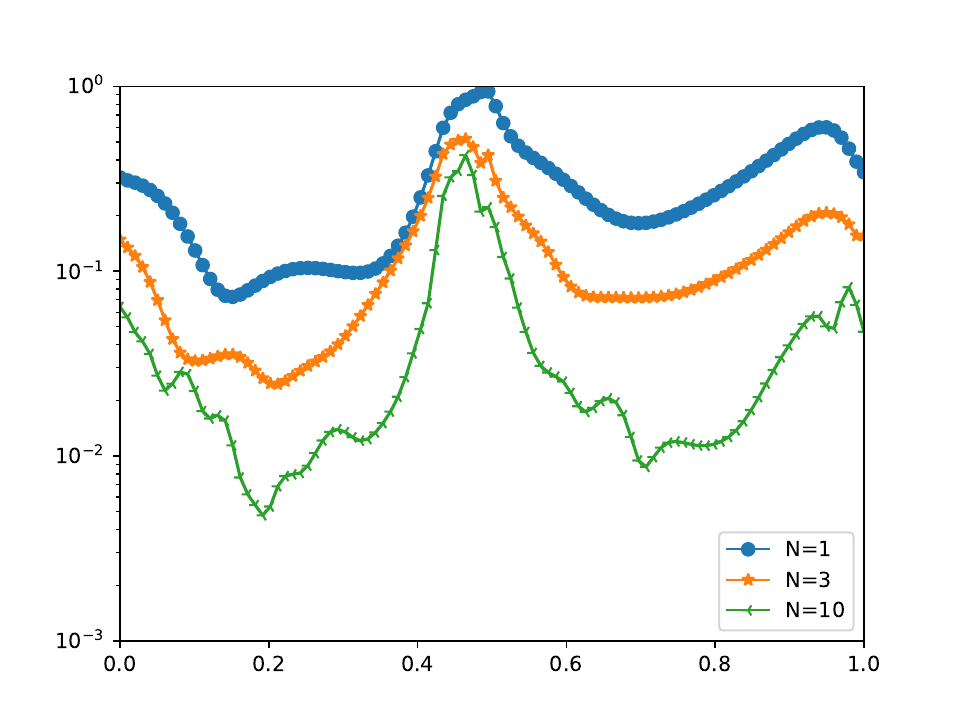}
    \put(47,0){Time}
    \put(-4,25){\rotatebox{90}{Relative error}}
    \end{overpic}}
 \caption{Cardiovascular benchmark: time evolution of the relative error for velocity (left) and WSS (right) at increasing of the number of modes.}
 \label{fig:errors_cabg}
 \end{figure}
\begin{figure}
	\centering
    \subfloat[][FOM \label{fig:ana-fom}]{
    \begin{overpic}[height=0.4\textwidth]{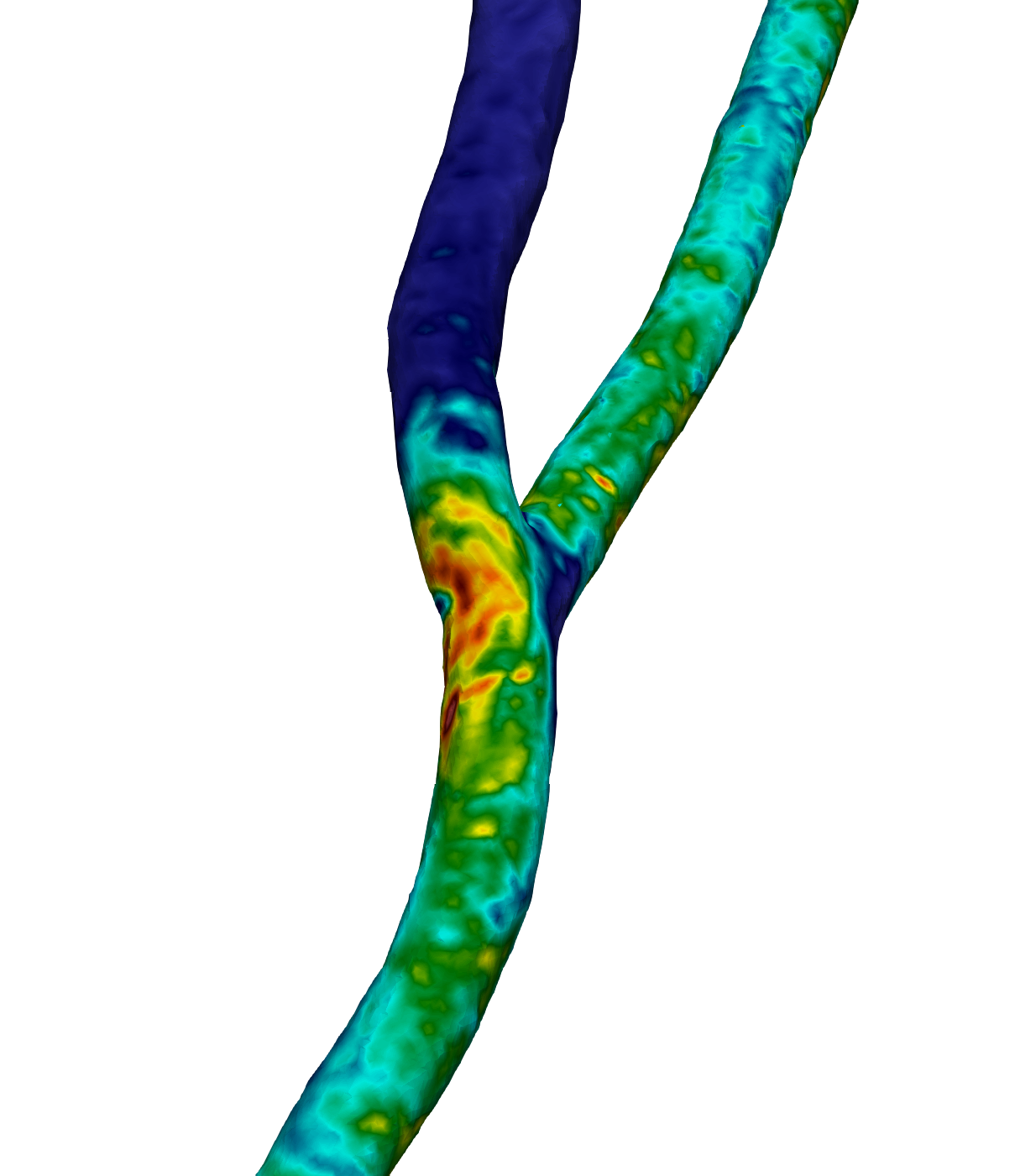}
    \end{overpic}}
    \subfloat[][ROM \label{fig:ana-rom}]{
    \begin{overpic}[height=0.4\textwidth]{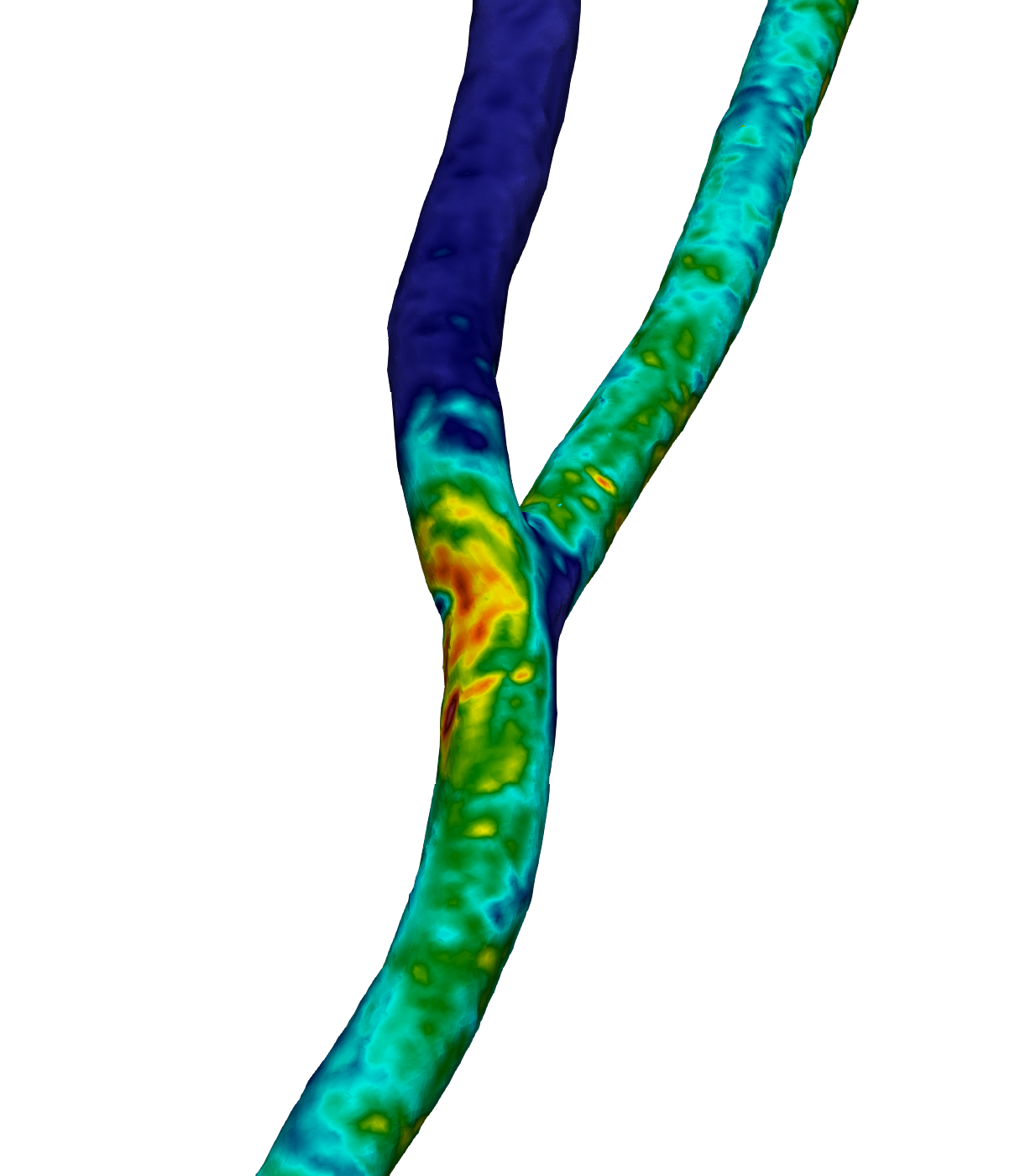}
    \end{overpic}}
    \hspace{2ex}
    \begin{overpic}[height=0.25\textwidth]{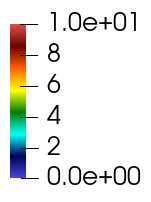}
    \put(-10,25){\rotatebox{90}{WSS (Pa)}}
    \end{overpic}
    

 \caption{Cardiovascular benchmark: qualitative comparison between WSS  distribution in the anastomosis
region computed by the FOM and by the ROM for the test point
(70\% stenosis) at $t = 0.8$ s.}
 \label{fig:qual-anast}
 \end{figure}
\begin{figure}
	\centering
    \subfloat[][FOM \label{fig:sten-fom}]{
    \begin{overpic}[height=0.4\textwidth]{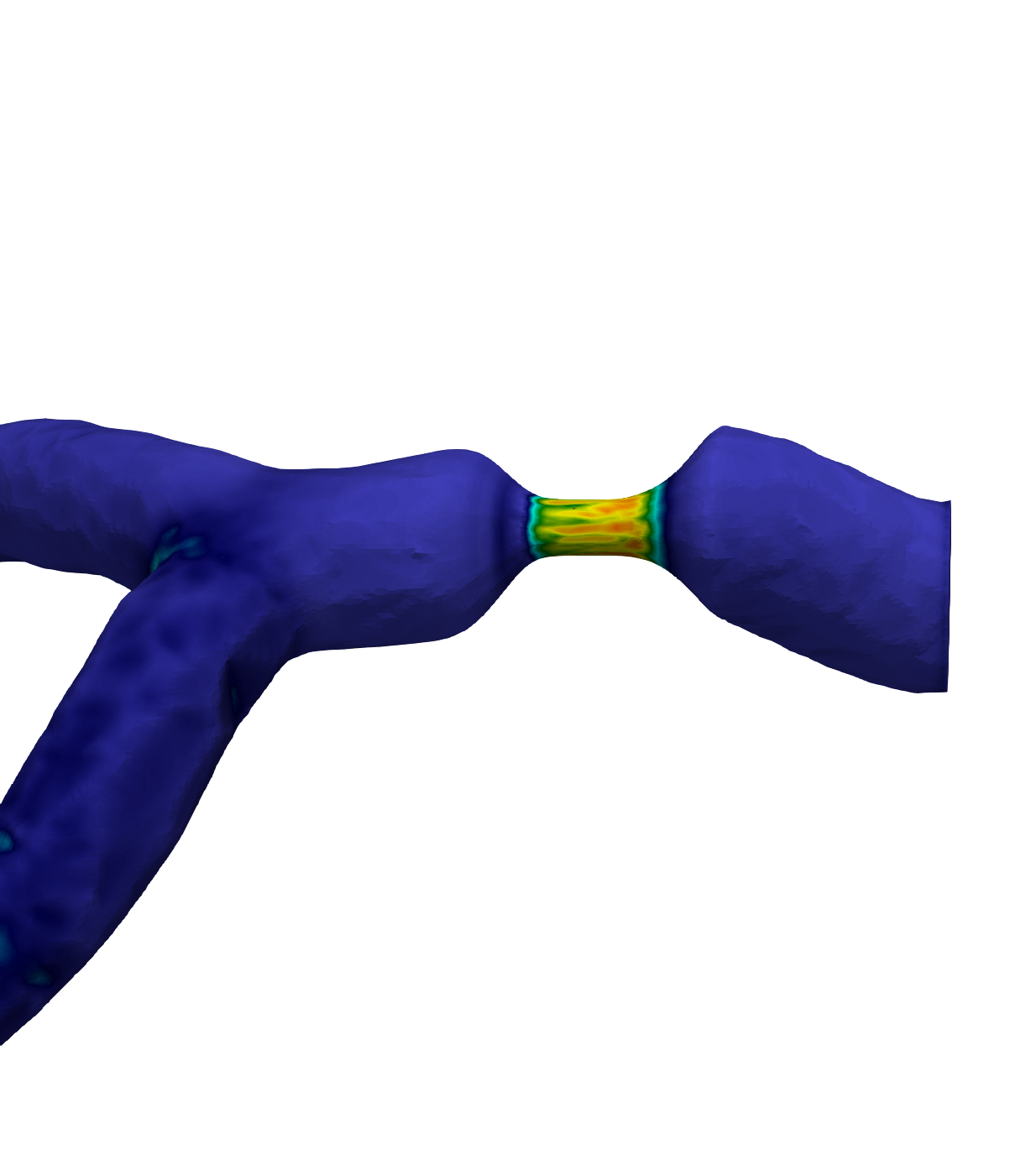}
    \end{overpic}}
    \subfloat[][ROM \label{fig:sten-rom}]{
    \begin{overpic}[height=0.4\textwidth]{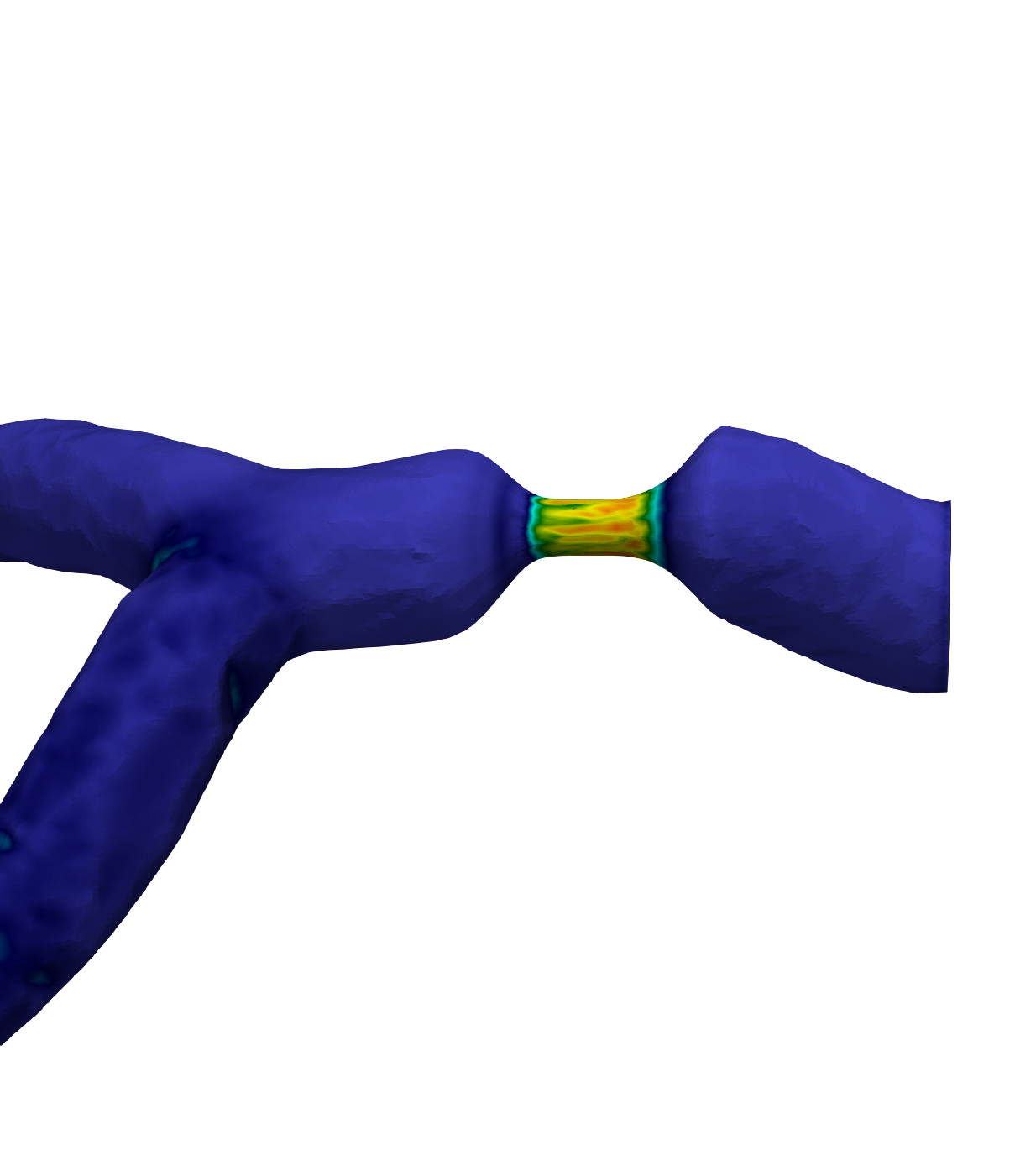}
    \end{overpic}}
    \hspace{2ex}
    \begin{overpic}[height=0.25\textwidth]{img/legenda.png}
    \put(-10,25){\rotatebox{90}{WSS (Pa)}}
    \end{overpic}
    

 \caption{Cardiovascular benchmark: qualitative comparison between WSS distribution in the stenosis
region computed by the FOM and by the ROM for the test point
(70\% stenosis) at $t = 0.8$ s.}
 \label{fig:qual-sten}
 \end{figure}

Trial and error process is employed to optimize the hyperparameters of the neural network used to map the reduced coefficients. The optimized architecture consists of 3 hidden layers, 1300 neurons per layer, a learning rate of $10^{-6}$ and hyperbolic tangent (Tanh) as activation function. 
After $50.000$ epochs, we reach an accuracy of about $93\%$. The ROM is trained based on a uniform sample distribution of the stenosis severity ranging between $50\%$ to $75\%$ with step $5\%$, except $70\%$ which is considered as validation point. 100 snapshots in time are collected for each value of stenosis severity, which results in a total of 500 training snapshots. 
Figure~\ref{fig:errors_cabg} shows the time evolution of the relative error for the velocity and WSS at varying of the number of POD modes. As expected, the relative error monotonically decreases at increasing of the number of modes. A time-averaged error of about $3.8\%$ and $4.9\%$ is obtained for the velocity and WSS, respectively, using 10 modes. \textcolor{black}{
} 
\textcolor{black}{Concerning an illustrative comparison, we limit to consider WSS due to its association with the restenosis process related to the stenosis and anastomosis regions. 
From Figures \ref{fig:qual-anast} and \ref{fig:qual-sten} we can observe that 
FOM and ROM solutions are very similar.  
}

Now we  briefly discuss the computational efficiency of our data-driven ROM. A FOM simulation takes about 7 hours. On the other hand, the computation of the reduced coefficients takes about 10 s. Then the resulting speed-up is about $10^5$ which a very relevant value, significantly greater than the one obtained with an intrusive ROM applied to a toy model (see Sec. \ref{sec:termal}). This proves the capability of data-driven ROMs to achieve real-time predictions in real-world applications.


\subsection{An industrial application}\label{sec:industry}
Now we  present an example coming from pharmaceutics industry in the framework of a research collaboration with Dompé farmaceutici S.p.A. concerning a system of gas-solid fluidized bed where the small particles are transported by the carrier gas flow \cite{goldschmidt2001hydrodynamic, fernandes}. Such a problem involves a multiphysics scenario and it is useful to show the performance in terms of accuracy and efficiency of data-driven ROM approaches for complex cases. See also Sec. \ref{sec:argos}.

The mathematical modelling is based on the so-called Computational Fluid Dynamics-Discrete Element Method (CFD-DEM) technique where the fluid and solid phases are described based on Eulerian and Lagrangian approach, respectively. The coupling between phases is given by the fluid-particle interaction forces \cite{Xu, tsu}. 

The FOM in strong formulation is given by the volume-averaged Navier-Stokes equations which read \cite{clayton}:

\begin{equation}
\label{momentum}
\left\{    
\begin{alignedat}{3}
\frac{\partial \epsilon}{\partial t} + \nabla \cdot ( \epsilon \mathbf{u} ) & = 0, & \quad & \text{in} \ \Omega \ \times \ (t_0, T],
\\
\frac{\partial (\epsilon \mathbf{u})}{\partial t} + \nabla \cdot ( \epsilon \mathbf{u} \otimes \mathbf{u}) & = - \nabla p - S_p + \nabla \cdot (\epsilon \bm{\tau}) + \epsilon \textbf{g}, & \quad & \text{in} \ \Omega \ \times  \ (t_0, T].
\end{alignedat}
\right.
\end{equation}

In the system above, $\epsilon$ is the fluid volume fraction and 
$\bm g$ is the gravity. 
The viscous stress tensor $\bm{\tau}$ is computed as follows: 

\begin{equation}
\bm{\tau} = \bm{\tau}_1 + \bm{\tau}_2 = \nu \left( \nabla\mathbf{u} + \nabla\mathbf{u}^T \right)- \dfrac{2}{3} \nu \left(\nabla \cdot \mathbf{u} \right) \bm{I},
\label{tau}
\end{equation}
where 
$\bm{I}$ is the identity matrix.

As for the cardiovascular application, a second-order finite volume method is adopted for the space discretization. The computational domain $\Omega$ is partitioned into $N_c$ cells or control volumes $\Omega_i$ with $i = 1, \dots, N_c$. 
The fluid volume fraction $\epsilon_{i}$, indicating the portion of the cell $i$ occupied by fluid, is defined as \cite{weller}: 

\begin{equation}
\epsilon_{i} = 1 - \frac {\sum_{j=1}^{n_p} \widetilde{\Omega}_j}{\Omega_{i}},
\label{eps}
\end{equation}
where  $n_p$ is the number of particles in the cell $i$ and $\widetilde{\Omega}_j$ represents the volume of the particle $j$. The interaction between fluid and particles is modelled by the source term $S_p$ which in the  cell $i$ is computed as \cite{zhu2007discrete}: 

\begin{equation}
S_{p,i} = \frac{\sum_{j=1}^{n_p} ( \mathbf{F}_{d,j} + \mathbf{F}_{\nabla p,j})}{\rho_f{\Omega}_{i}}.
\label{source}
\end{equation}

Here, $\mathbf{F}_{d,j}$ and $\mathbf{F}_{\nabla p,j}$ are the drag and pressure gradient forces, respectively \cite{zhu2007discrete}. 




Like for the cardiovascular test, the time interval of interest $(t_0, T]$ is partitioned using the time step $\Delta t \in \mathbb{R}^+$, resulting in time levels $t^n = t_0 + n \Delta t$, where $n$ ranges from 0 to the total number of time steps $N_{T}$. The system~\eqref{momentum} is discretized in time employing a first-order Euler scheme. On the other hand, all terms in~\eqref{momentum} are treated implicitly except for the divergence term, which was discretized using a semi-implicit scheme. 


The DEM model 
resolves particle motion 
by the translation and rotation second Newton's laws which are described as follows \cite{review, molin, fernandes}:

\begin{equation}
m_j \frac{d\widetilde{\mathbf{u}}_j}{dt} =  \sum_{m=1}^{n_j^c} \bm{F}_{jm}^c +  \bm{F}_j^f + m_j \bm g ,\\
\label{translation}
\end{equation}
\begin{equation}
I_j \frac{d\mathbf \omega_j}{dt} = \sum_{m=1}^{n_j^c} \bm{M}_{jm}^c, 
\label{rotation}
\end{equation}
where 

\begin{equation}
m_j = \rho_p\dfrac{\pi d_p^3}{6} \quad \text{and} \quad I_j = \dfrac{m_j d_p^3}{6}
\end{equation}
are the mass and the moment of inertia of the  particle $j$, respectively. In addition, $\rho_p$,  $d_p$ and $\mathbf \omega_j$ are the density, the diameter and the angular velocity of the particle $j$, respectively.  $\bm{F}_{jm}^c$ and $\bm{M}_{jm}^c$ are the contact force and torque acting on particle $j$ by its $m$ contacts, which may involve other particles or wall \cite{fernandes}. The total number of contacts for the particle $j$ is denoted as $n_j^c$ and $\bm{F}_j^f = \bm{F}_{d,j} + \bm{F}_{\nabla p, j}$ is the particle-fluid interaction force acting on particle $j$. Notice that non-contact forces are not considered in this work.

The Stokes number, denoted as $Stk$, which characterizes the behavior of particles suspended in a fluid flow, is computed as: 
\begin{equation}
    Stk = \frac{\tau_p}{\tau_f},
\label{Stk}
\end{equation}
where $\tau_f$ is the carrier fluid characteristic time and $\tau_p$ is the particle relaxation time.
Eqs. \eqref{translation}-\eqref{rotation} are discretized by adopting a first-order Euler scheme.

For the multiphysics system at hand, we introduce a “local” variant of the “global” POD-RBF data driven approach described in Sec. \ref{sec:data-rom}. Unlike the method described in Sec. \ref{sec:data-rom} here a POD basis is computed for each parameter in the
training set and the basis functions for new parameter values are found via RBFs interpolation of the basis functions associated to the training set. 
It should be noted that local POD-RBF and global POD-RBF coincide when only time reconstruction is considered. For more details the reader is referred to \cite{hajisharifi2023non}.

The Stokes number ${Stk}$, defined in eq. \eqref{Stk}, represents a crucial parameter of the model. 
We build a parametric ROM with respect to $Stk$ for the Eulerian phase. 
We choose a uniform sample distribution in the range $Stk \in [200, 300]$. To vary $Stk$ we modify the value of $\rho_p$. We consider 21 sampling points. For each value of $Stk$ in such set, a simulation is run for the entire time interval of interest, i.e. $(0, 5]$. 
The snapshots are collected every 0.02 s, for a total number of 10500 snapshots.
To train the ROM we consider the 90$\%$ of the database for each simulation, i.e. 450 snapshots, randomly chosen,  resulting in a total of 9450 snapshots. The remaining ones are considered as validation set. 

\begin{figure}[ht]
\centering
\subfloat[ Global POD-RBF]{\label{fig:CumEnergy_global}\includegraphics[width=.52\linewidth]{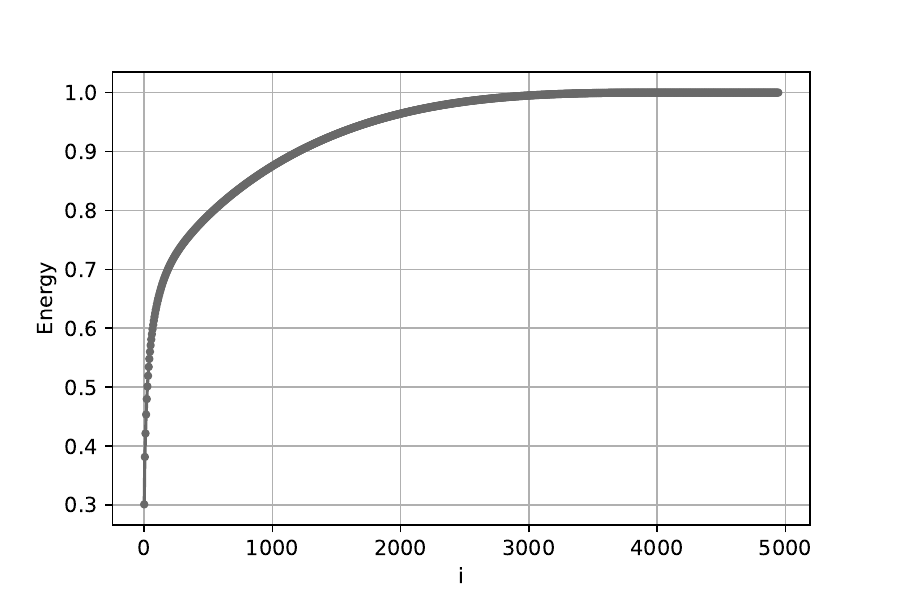}}
\subfloat[Local POD-RBF]{\label{fig:CumEnergy_local}\includegraphics[width=.52\linewidth]{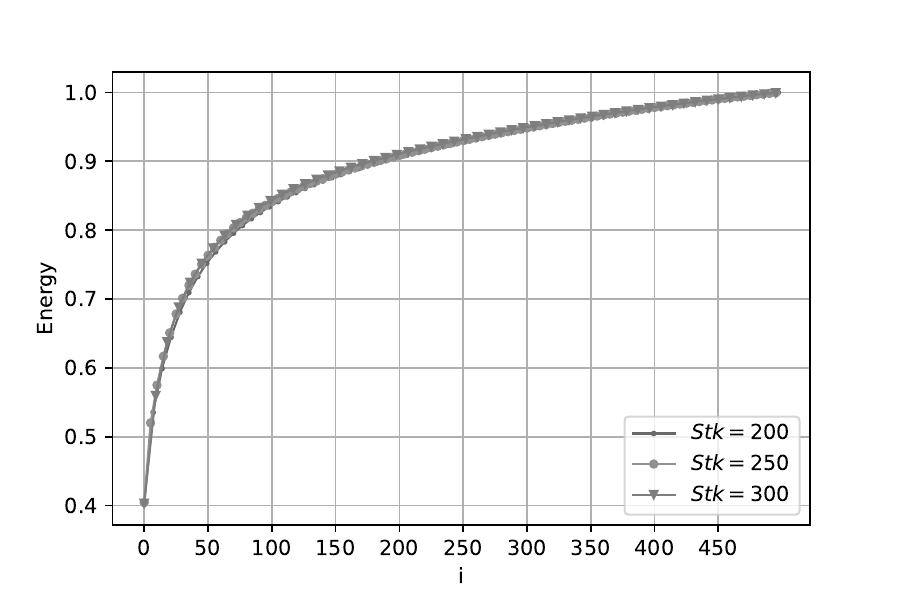}}
\caption{Industrial benchmark: cumulative eigenvalues for the fluid volume fraction $\epsilon$ by global POD-RBF (a) and local POD-RBF (b). }
\label{fig:Energy_local_Global}
\end{figure}

The plot of the cumulative eigenvalues for the fluid volume fraction $\epsilon$ is shown in Fig. \ref{fig:Energy_local_Global}. For the global POD-RBF approach (a), we have truncated the plot at 5000 modes, since the 99.99\% of the energy is already reached for 3200 modes. 
However, we oberve that at least 1500 modes are necessary to recover the 90\% of the cumulative energy. 
On the other hand, concerning the local POD-RBF (b), we show the plot of the cumulative eigenvalues for the initial, mid and final value of $Stk$ in the training set, 200, 250 and 300 respectively. 
The comparison shows no significant difference between the three curves suggesting that the same number of modes needs to be considered to capture the same energy threshold. In particular, 150 modes can capture the 90$\%$ of the energy for the whole range of $Stk$, one order of magnitude in less than the global POD-RBF.

\begin{figure}
    \centering
    \subfloat[POD-RBF validation at $t =1$ s.]{\label{fig:St227_t=1}\includegraphics[width=10cm]{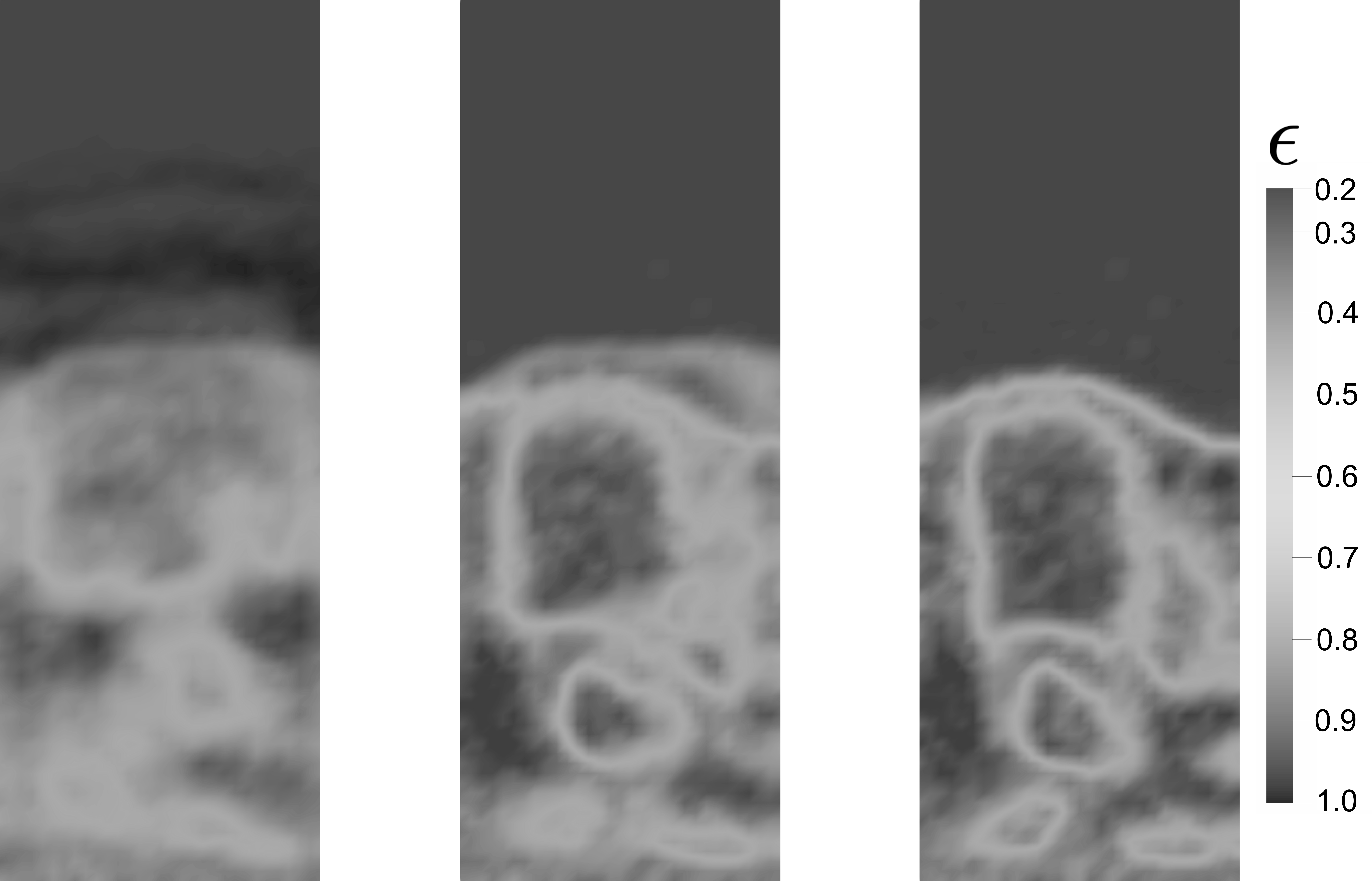}}\\
    \subfloat[POD-RBF validation at $t =2.5$ s.]{\label{fig:St227_t=2.5}\includegraphics[width=10cm]{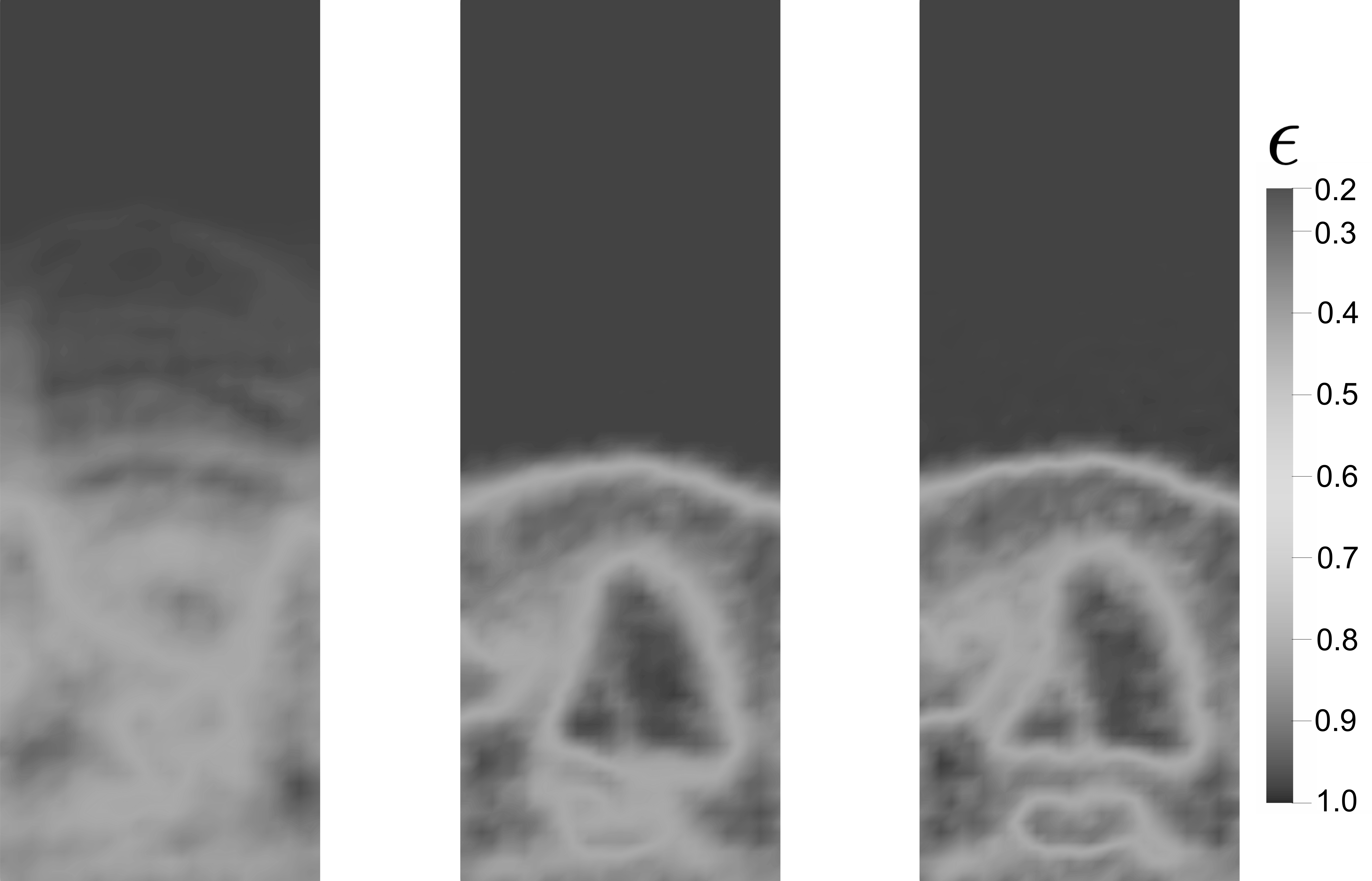}}\\
   \subfloat[POD-RBF validation at $t =4$ s.]{\label{fig:St227_t=2.4}\includegraphics[width=10cm]{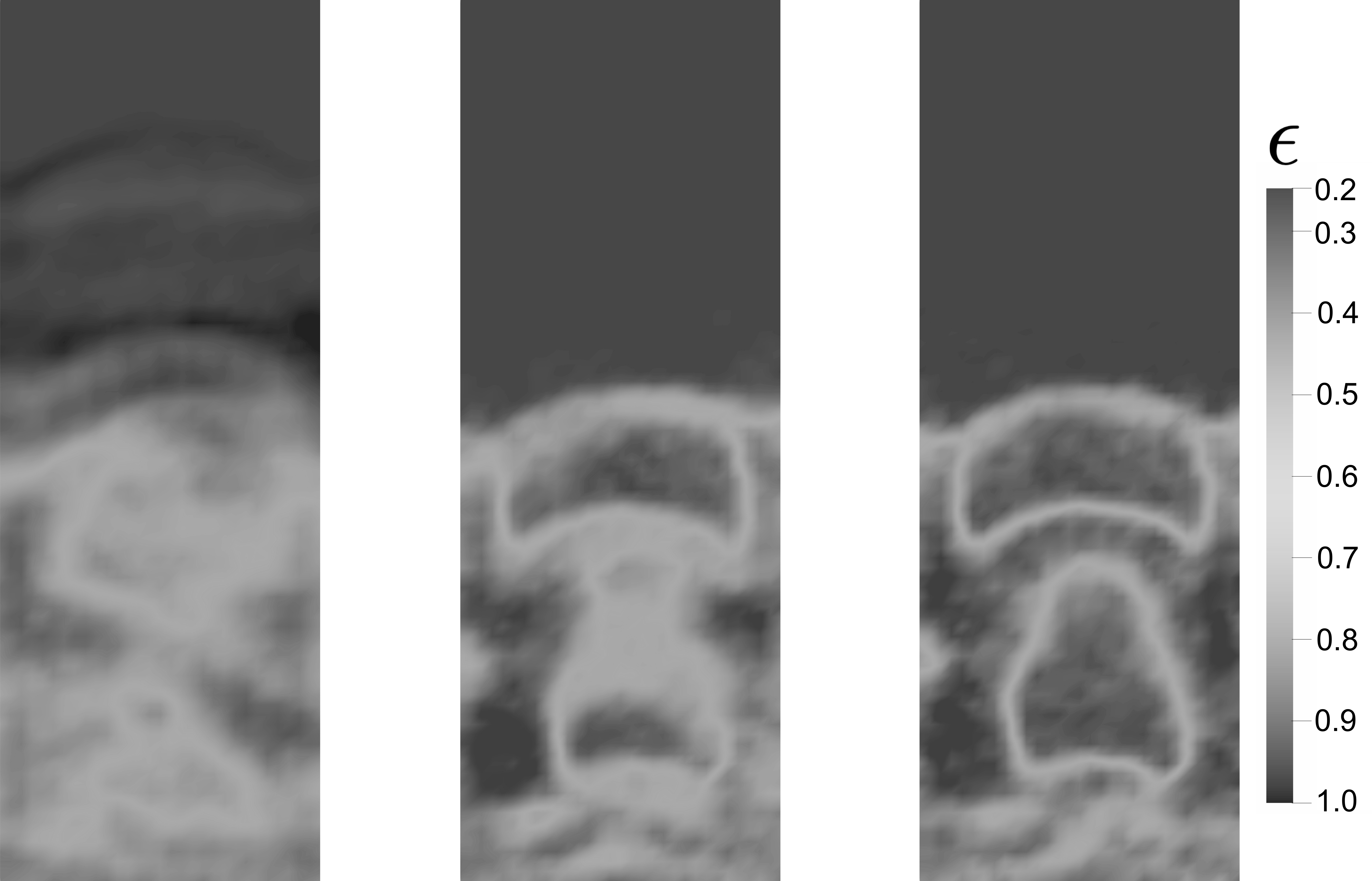}}\\
    \caption{Industrial benchmark: comparison between FOM and ROM for $Stk =227$ for three different time instances: $t = 1$ s (first raw), $t = 2.5$ s (second raw) and $t = 4$ s (third raw). The first column represents the solution by global POD-RBF, the second one the solution by local POD-RBF and the third one the solution by FOM.} 
    \label{fig:Vis_ParROM_St227}
\end{figure}

We take $Stk = 227$ as testing point to
evaluate the performance of the ROM. A comparison between the two reconstructed solutions and the corresponding FOM is shown in Fig. \ref{fig:Vis_ParROM_St227}  for $t = 1$ s, $t = 2.5$ s and $t = 4$ s. 
As one can see, in both cases, the global POD-RBF is barely able to reconstruct the main patterns of flow field. The solution is affected by some spurious oscillations. The situation changes using the local POD-RBF which can capture more details and partially damps the unphysical oscillations. 

\begin{figure}[ht]
\centering
\subfloat[ $Stk$ = 227]{\label{fig:Stk227_error}\includegraphics[width=.62\linewidth]{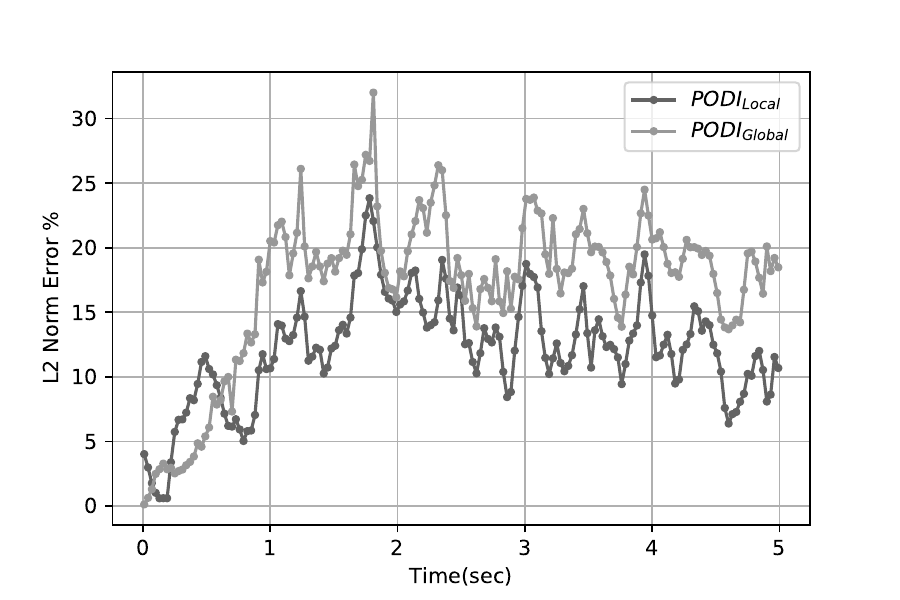}}
\caption{Industrial benchmark: time evolution of the relative error for the testing point $Stk$ = 227. In the acronym PODI, “I” means “Interpolation” performed by RBFs.}
\label{fig:error_time_parametric}
\end{figure}

To introduce a more quantitative comparison, 
the time history of the  relative error between FOM and ROM has been computed for $Stk = 227$. The plot is shown in Fig. \ref{fig:error_time_parametric} whilst Table \ref{tab:Parametric_error} reports the values of the maximum, minimum and mean error. From Fig. \ref{fig:error_time_parametric}, we note that, except for the initial time of the simulation, the local POD-RBF provides an error lower than the global one. The maximum error decreases of about two times when one adopts the local POD-RBF. Also the mean error decreases passing from 17-18\% to 12\%. 




Beyond the capability of local POD-RBF to capture better the system dynamics, it is also much cheaper in terms of CPU time of about one order of magnitude.  The FOM
simulation takes around 1.8e5 s. The online phase associated the local POD-RBF takes around 9 s, while the global POD-RBF requires 121 s. So the computational speed-up is of order 1e4 for local POD-RBF and 1e3 for the global one. 

So we can conclude that the local POD-RBF basically performs better than the global one both in terms of efficiency and accuracy. This result confirms the goodness of data-driven surrogate models for complex applications.


\begin{table}[ht!]
\centering
\begin{tabular}{||c|c|c||}
\hline 
$Stk$ &  \multicolumn{2}{c|}{$227$}  \\
\hline\hline
ROM approach &  Local POD-RBF & Global POD-RBF \\
\hline
Mean Error (\%) &  12.2 & 18.1  \\ \hline
Max Error (\%) & 19 & 33 \\\hline
Min Error (\%) &4.3 & 0.6  \\\hline
\end{tabular}
 \caption{Industrial benchmark: mean, maximum and minimum values of the relative error for the testing point $Stk = 227$ by local and global POD-RBF.}
 \label{tab:Parametric_error}
\end{table}




\section{Argos and Atlas}
\label{sec:argos_atlas}

Over the past few decades, there's been a remarkable boost in numerical simulations, thanks to advancements in numerical algorithms, computational mathematics, and hardware technology. Despite the progress with increasingly sophisticated algorithms and more powerful computers, many challenges still persist. Complex problems still demand considerable time and computational resources. This is where ROMs step in. As illustrated in the previous sections, such techniques have become quite popular across engineering and applied science fields for their ability to quickly approximate solutions to Parametric Partial Differential Equations in real-time. However, the lack of easy-to-use, readily available (even commercial) applications for conducting ROM simulations is holding back a wider adoption.

While implementing ROMs introduces complexities in ensuring their fidelity and reliability, recent advancements have demonstrated notable successes in leveraging ROMs for real-world simulations, particularly in industrial and cardiovascular domains.

During the ERC CoG AROMA-CFD project \cite{aromabook}, several certified, precise, and rapid methodologies have been developed and implemented. The goal was to create an efficient computational framework capable of analyzing parametric flow problems in real-time, overcoming various methodological barriers and limitations. However, these methods are currently housed in generic libraries, necessitating programming skills and a significant time investment to configure the program for specific problems.

To bridge this gap, we have developed two Proof of Concept web platforms, namely Argos (\url{https://argos.sissa.it/}) and Atlas (\url{https://atlas.sissa.it/}), that aim to meet the demand for ready-to-use tools for industrial and cardiovascular real-time simulations, respectively. Notably, the accessibility of ROM-based simulations through web platforms has further democratized their utilization, allowing practitioners to harness the benefits of reduced order modeling without the burden of intricate implementation details.

Through a critical examination of these advancements and their implications, we aim to provide insights into the latest technological advancements of ROMs, highlighting both their potential and the remaining hurdles on the path towards their widespread adoption and seamless integration into computational workflows. By doing so, we seek to elucidate the intricate interplay between theoretical developments, technological innovations, and practical applications in the ever-evolving landscape of computational modeling and simulation.

\subsection{Argos}
\label{sec:argos}
A generic simulation aimed at optimizing a particular industrial object or artifact typically involves several complex phases, such as digitally generating the object's geometry and subsequently evaluating it using a traditional numerical solver. However, achieving accurate results often requires multiple iterations, consuming significant computational and human resources. Moreover, configuring the numerical framework and managing pre- and post-processing phases further prolong the analysis time.

\begin{figure}[ht]
\centering
\includegraphics[width=.95\linewidth]{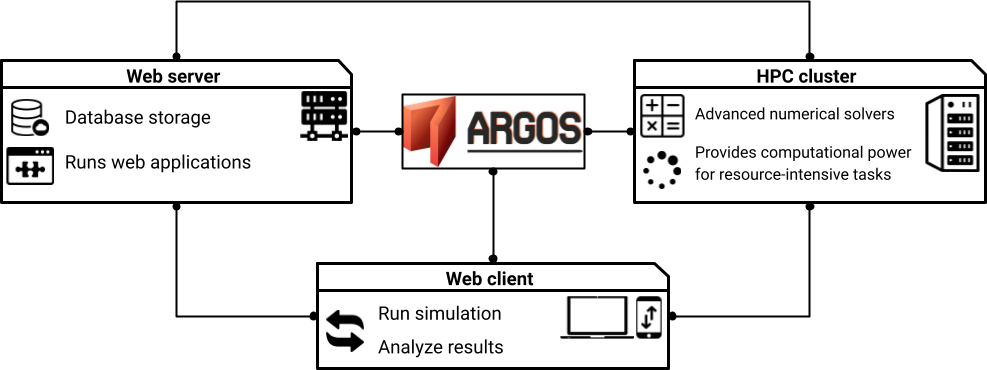}
\caption{Argos infrastructure: interplay between the web server, the backend High-Performance Computing (HPC) cluster, and the web client.}
\label{fig:argos}
\end{figure}

\begin{figure}[ht]
\centering
\includegraphics[width=.95\linewidth]{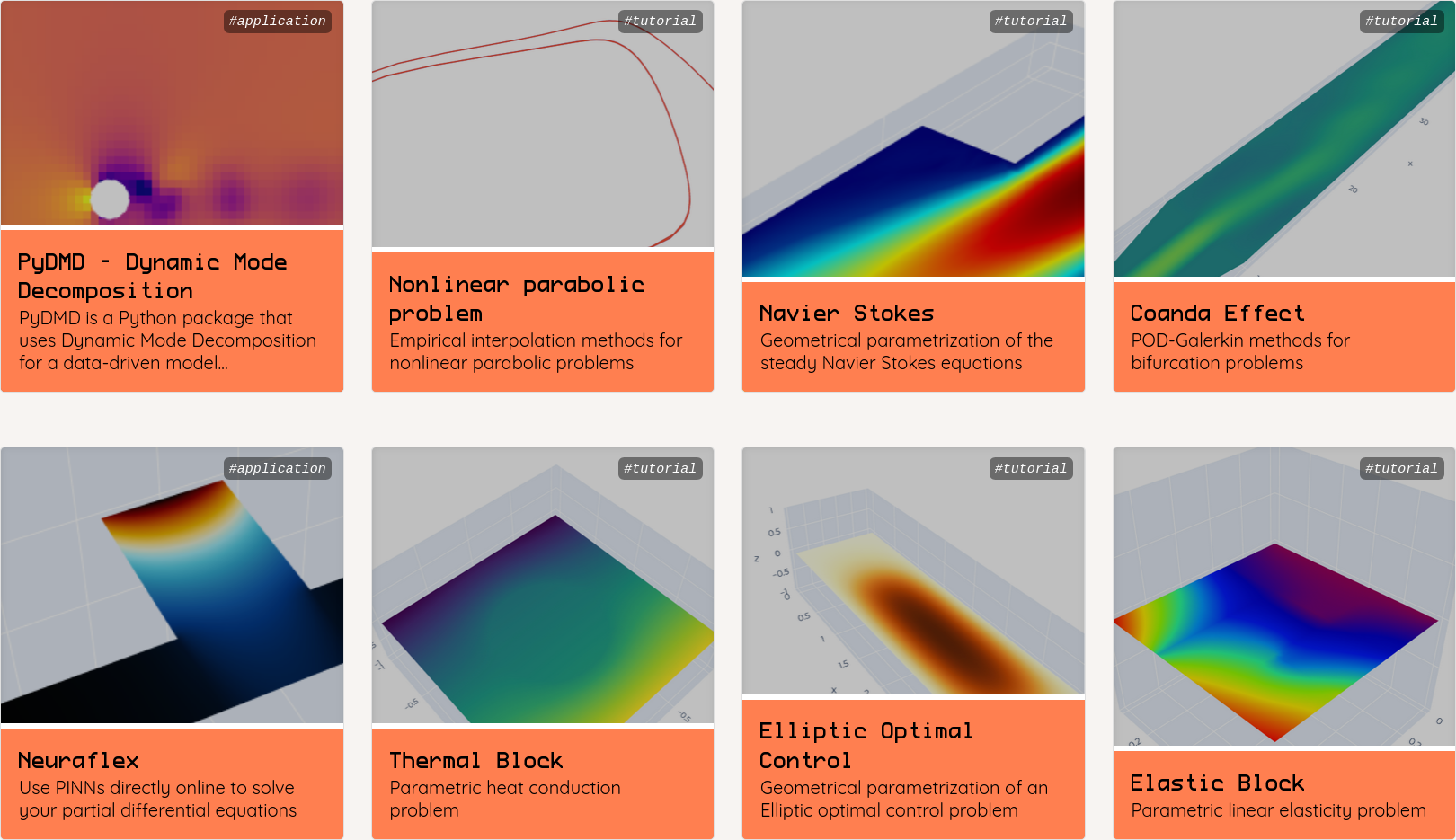}
\caption{Selection from the gallery of applications and tutorials available on the Argos web platform.}
\label{fig:argos_gallery}
\end{figure}

\begin{figure}[ht]
\centering
\includegraphics[width=.95\linewidth]{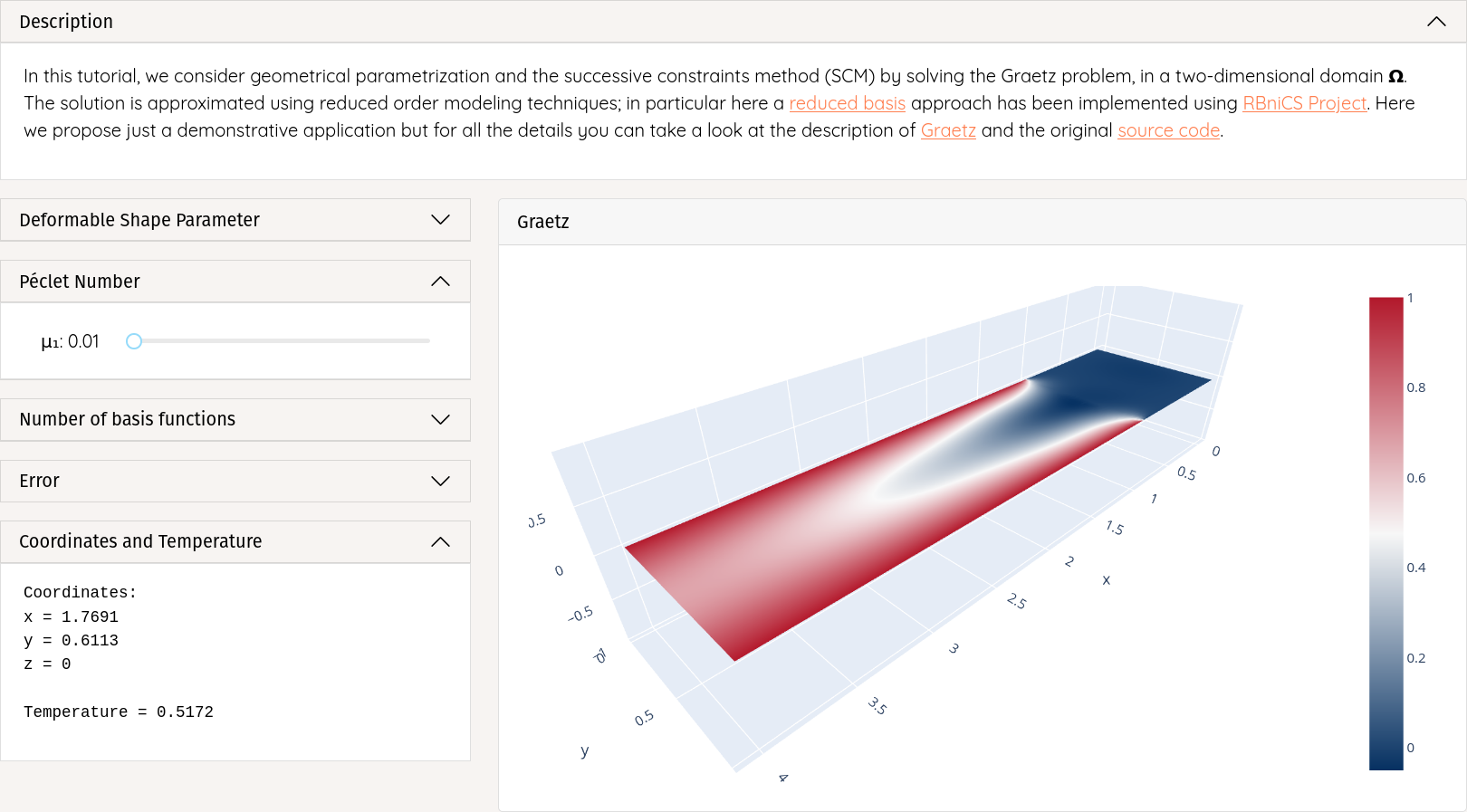}
\caption{Example of a typical graphical user interface for an Argos application. This particular example demonstrates the solution to the Graetz problem. On the left side of the interface, users can manipulate various parameters, including the Péclet number and the number of basis functions utilized for the ROM expansion. They can also visualize the solution field at a specified point. On the right side, the solution across the entire domain is displayed in an interactive plot.}
\label{fig:argos_example}
\end{figure}

Argos (\url{https://argos.sissa.it/}) aims to revolutionize simulation tools for industrial applications, offering a new paradigm where complex parametric analyses are portable, intuitive, and real-time. With Argos, simulations are accessible via the internet, as illustrated in Fig.~\ref{fig:argos}. Users can simply select the desired application from the provided showcases, such as those shown in Fig.~\ref{fig:argos_gallery}, and fine-tune it with the appropriate physical quantities describing the object or artifact. A concrete example is displayed in Fig.~\ref{fig:argos_example}. Overall, the long waits for simulation results and the need for extra activities are replaced by a more sustainable and efficient workflow.

While commercial solutions with similar features exist, as far as we know, they aren't consolidated into a single web-accessible portal. Additionally, many existing solutions don't leverage state-of-the-art methodologies. Indeed, the backbone of Argos server is built upon the libraries developed during the ERC AROMA-CFD project, representing some of the most advanced frameworks for reduced order modeling.

Leveraging state-of-the-art developments in data-driven reduced order models, Argos encapsulates these methodologies into user-friendly graphical applications. The goal of Argos is to democratize numerical simulation, making it accessible to a wide audience of companies and institutions through simple web interfaces. By providing graphical and ready-to-use applications involving ROM, we addressed a huge gap in the market.

The web interface of Argos ensures maximum portability across various hardware and operating system environments. This approach allows for the development of a single version for each application, regardless of the digital ecosystems of potential customers. This streamlines the development phase, making it faster and more sustainable, while enhancing portability through internet protocols accessible via standard web browsers.

Employing ROM techniques also creates synergies with the client-server architecture. The low computational cost inherent in ROM allows applications to run efficiently on a single server without encountering computational bottlenecks. The algorithmic computations of ROM are handled server-side, significantly reducing hardware requirements on the client side and promoting a more efficient and sustainable framework. Additionally, server-side management of installations, upgrades, and configurations of required libraries by administrators eliminates the need for end-users to manage these tasks.

For programming language, Python was chosen to maximize compatibility with existing software. This allows for seamless integration of AROMA-CFD numerical tools, which typically expose Python Application Programming Interfaces (APIs), with a graphical framework. This approach eliminates the need to bind multiple programming languages, significantly reducing human resources required for application development. Leveraging Python also provides access to several advanced computational libraries for scientific computing, numerical analysis and data science. These technical choices enable Argos to simplify and expedite the development of new applications, ensuring competitiveness in the market despite differences in human resources between the Argos team and other development-focused companies.

Argos has become a pivotal computational tool in facilitating the integration of High-Performance Computing (HPC), reduced order modeling, data science, Internet of Things (IoT), artificial intelligence, and digital twin developments. Its remarkable versatility and innovative potential are firmly grounded in robust methodological evidence, positioning it as a cornerstone in a landscape where models, data, and digital technology are radically transforming computational science and engineering. Furthermore, the Argos infrastructure holds the potential for extensive future utilization, providing a platform to host cutting-edge developments in simulation. This includes emerging approaches such as physics-informed models and quantum computing, ensuring Argos remains at the forefront of computational innovation in the long term.

\subsection{Atlas}
\label{sec:atlas}
The human cardiovascular system is an incredible network of blood vessels, chambers, valves, mechanical deformations, and electrical signals. Together, they work seamlessly to pump oxygenated blood throughout the body. But with this intricate design comes vulnerability. The cardiovascular system is susceptible to a range of diseases and conditions, including coronary artery disease, heart failure, arrhythmias, and hypertension. These conditions are among the leading causes of mortality and morbidity worldwide.

In the context of cardiovascular research and clinical practice, one of the biggest challenges is understanding and predicting how this complex system behaves under various physiological and pathological conditions, especially in individual patients. This is where CFD simulations come into play. They're a powerful tool for modeling and analyzing the intricate flow patterns, pressure gradients, and shear stresses that impact vascular health. By providing such accurate information, CFD simulations offer valuable insights that can drive medical advancements, refine treatment strategies, and deepen our understanding of this vital physiological process.

However, despite their potential, current cardiovascular simulations face several limitations, some of which already mentioned in Sec.~\ref{sec:coronary}. Firstly, high-fidelity simulations often require significant computational resources and time, making them impractical for real-time clinical decision-making. Then there's the issue of personalization—generic models don't consider individual patient variations, which limits the ability to tailor treatments and interventions. Additionally, while highly complex simulations may yield accurate results, they often lack interpretability, making it difficult for healthcare professionals to extract meaningful insights and practical guidelines. Lastly, these cutting-edge simulations are often confined to specialized research environments, restricting access to a wider audience.

\begin{figure}[ht]
\centering
\includegraphics[width=.95\linewidth]{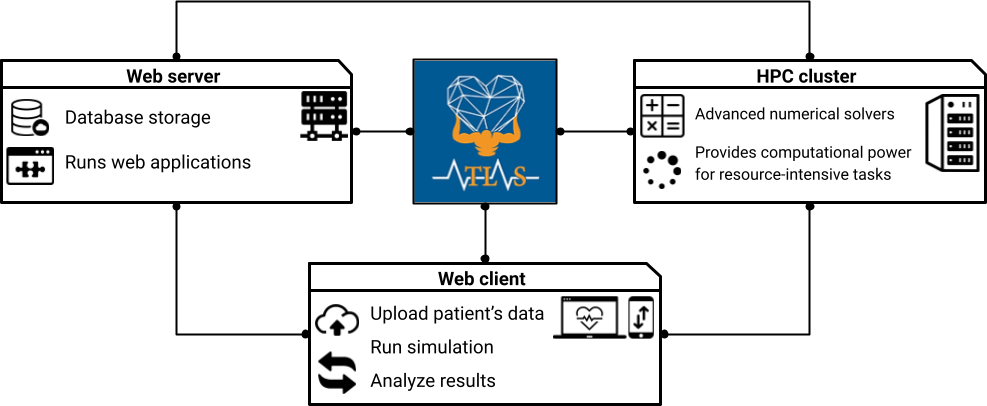}
\caption{Atlas infrastructure: interplay between the web server, the backend High-Performance Computing (HPC) cluster, and the web client.}
\label{fig:atlas}
\end{figure}

In response to the challenges and limitations faced in cardiovascular simulations, we have created Atlas (\url{https://atlas.sissa.it/}), a groundbreaking cloud computing platform. Atlas aims to revolutionize our approach to cardiovascular simulations by providing clinicians, researchers, and healthcare professionals with access to highly accurate and easily interpretable live simulations. The platform will be directly accessible via web browsers, eliminating the need for deep technical expertise in mathematical modeling or programming. Atlas is designed to be highly accessible. As a web-based platform (Fig.~\ref{fig:atlas}), it will be easily accessible to a wide range of users, including clinicians, researchers, and healthcare professionals. This accessibility fosters collaboration and knowledge dissemination across different disciplines, ultimately contributing to advancements in cardiovascular modeling and healthcare decision-making.

The development of Atlas is driven by the need to bridge the gap between complex simulations and real-world clinical requirements. One of its primary objectives is to enhance efficiency by leveraging ROM techniques, which have been refined through previous research. By doing so, Atlas will be capable of delivering high-fidelity simulations with significantly reduced computational demands, enabling real-time access to critical analyses. This efficiency is further supported by a robust software infrastructure developed as part of the ERC AROMA-CFD project.

\begin{figure}[ht]
\centering
\includegraphics[width=.95\linewidth]{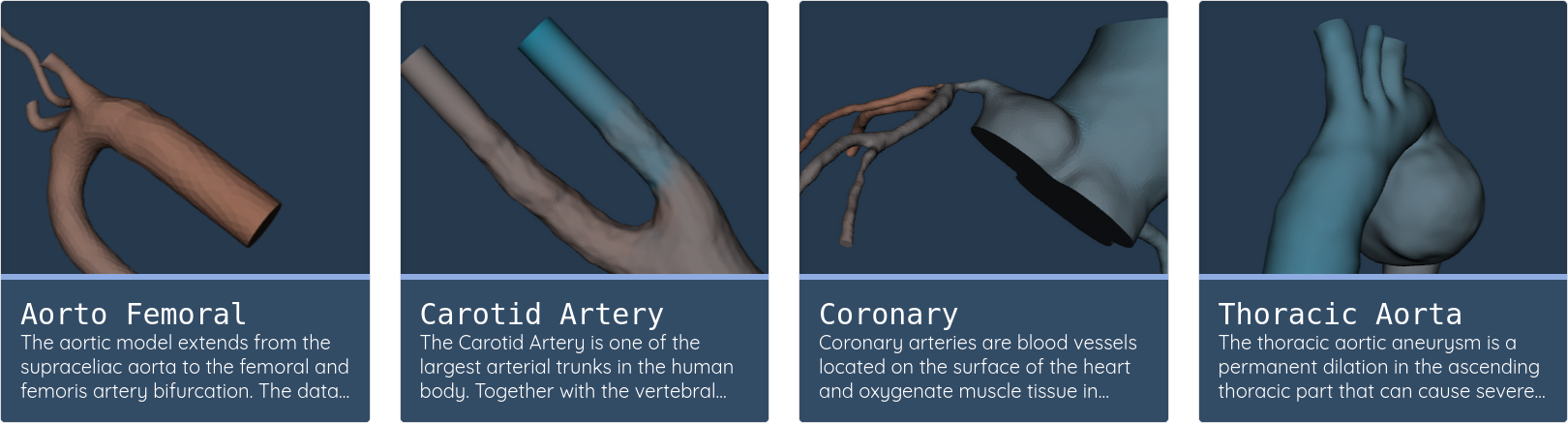}
\caption{Selection from the gallery of applications available on the Atlas web platform.}
\label{fig:atlas_gallery}
\end{figure}

As for Argos, users can simply select a desired application from the provided showcases shown in Fig.~\ref{fig:atlas_gallery}, and fine-tune it with the appropriate physical and patient-specific information. Another critical aspect of Atlas is its focus on personalization. Through the integration of patient-specific data, such as medical imaging and clinical measurements, the platform tailors simulations to individual profiles. This personalized approach, facilitated by data-driven techniques, empowers healthcare professionals to make informed decisions and develop personalized treatment strategies based on the unique characteristics of each patient, streamlining workflows and embracing the digital transformation of the healthcare sector.

\begin{figure}[ht]
\centering
\includegraphics[width=.95\linewidth]{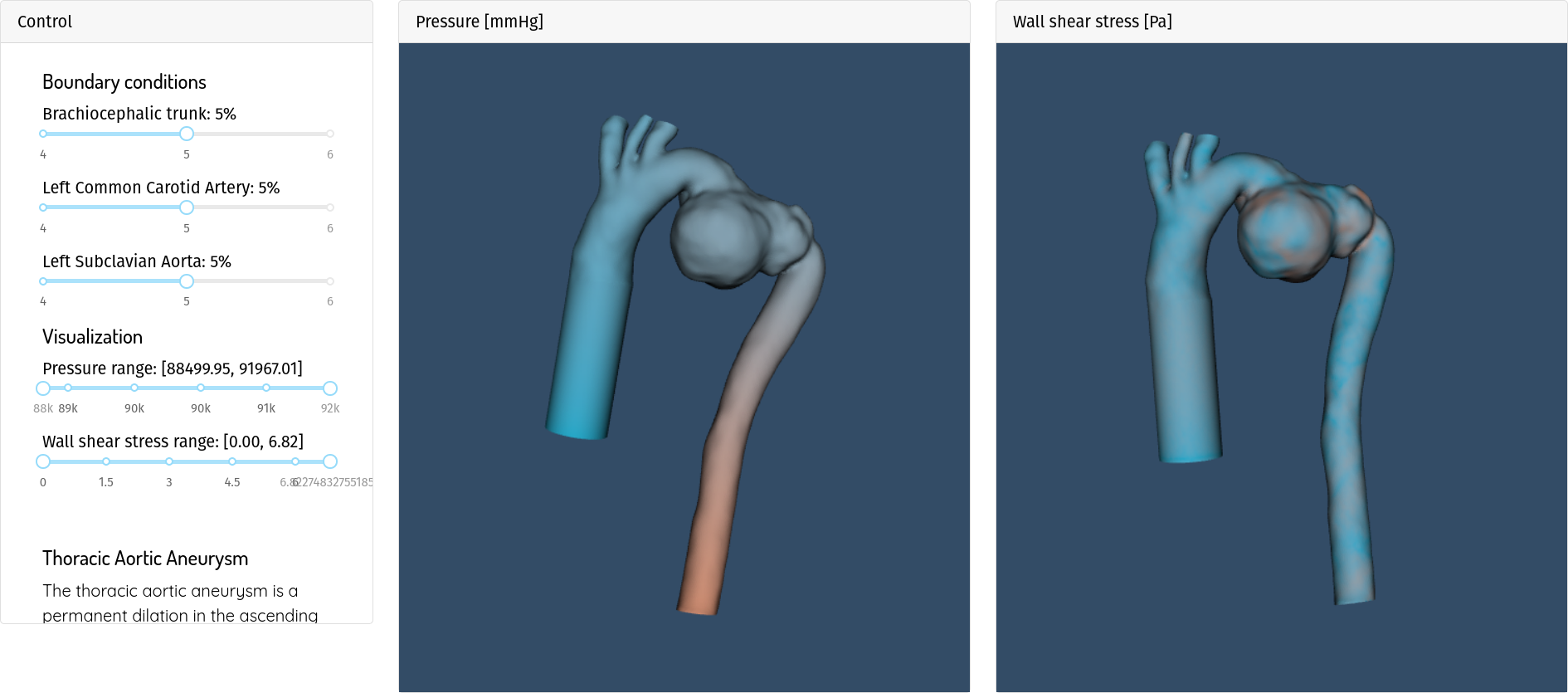}
\caption{Example of a typical graphical user interface for an Atlas application. This particular example demonstrates the blood flow field in the Thoracic Aorta of a patient suffering from aneurysm.. On the left side of the interface, users can manipulate various parameters, including the patient-specific boundary conditions at the different branches, typically derived from clinical measurements. On the right side, the pressure field and the Wall Shear Stresses (WSS), a predictor for aneurysm risk of rupture, across the entire domain are displayed in an interactive plot.}
\label{fig:atlas_example}
\end{figure}

Moreover, Atlas places a strong emphasis on interpretability. The platform is designed to ensure that the models and outcome analyses it provides are easily understandable to healthcare professionals. This transparency and interpretability are essential for building trust in the insights derived from the simulations, ultimately enhancing their utility in clinical decision-making. An example showcasing the blood flow in a thoracic aortic aneurysm is displayed in Fig.~\ref{fig:atlas_example}.

Existing solutions in the market, including commercial software, have been utilized for cardiovascular simulations. These libraries offer similar capabilities for advanced use cases in cardiovascular modeling. However, none of these solutions seem to be integrated into a single web-accessible portal for real-time, personalized model simulations. Additionally, they often lack reliance on state-of-the-art cloud computing and IoT methodologies. Notably, Atlas is built upon libraries developed during the ERC AROMA-CFD project, representing some of the latest advancements in ROMs.

The project outcomes of Atlas are innovative and distinctive due to several key factors. Atlas combines ROM with data-driven methods, offering a unique fusion that delivers the accuracy of complex simulations while minimizing computational demands. It provides real-time accessibility to high-fidelity simulations via a web browser, setting it apart from traditional tools that require substantial computational resources and time. Atlas underscores its ability to tailor simulations according to individual patient profiles, a level of personalization frequently absent in current solutions, rendering it indispensable.

Envisioning a futuristic scenario, integrating Atlas with Augmented Reality (AR) in surgery rooms opens up exciting possibilities. Surgeons could wear AR headsets that provide real-time patient-specific data and simulations from Atlas directly in their field of view during surgery. This futuristic synergy among ROMs, Atlas and AR holds the potential to revolutionize surgical practice and patient care, offering enhanced precision, decision support, personalized planning, and improved training and education.

\section{Conclusions}

In this chapter we have explored the accuracy and efficiency of ROMs in the context of academic and real-world applications, covering several techiques. 
In particular, we have introduced classical intrusive approaches, highlighting their feature to build surrogate models based on physical principles. Subsequently, we have explored data-driven models, which leverages machine learning and regression techniques to build ROMs directly from data. We have also presented some techniques for the computation of the reduced space, such as POD and greedy algorithms, which play a crucial role both for intrusive and non intrusive approaches.

Spanning various scenarios and test cases, we evaluated the performance of ROMs in different settings, ranging from an academic benchmark, to biomedical and industrial applications. In Sec. \ref{sec:termal}, the thermal block problem is treated by using an intrusive framework 
offering valuable insights into the accuracy and efficiency of this approach. \textcolor{black}{In this scenario, achieving errors smaller than $10^{-6}$ is feasible due to the simplicity of the benchmark.}
Then, the performance of data-driven ROMs have been analyzed for a patient-specific cardiovascular application focusing on coronary arteries in a geometrical parameterization context (Sec. \ref{sec:coronary}) and for 
an industrial problem related to the granulation process in pharmaceutical applications involving a fluid-solid system (Sec. \ref{sec:industry}). 
\textcolor{black}{Both of these cases exhibited greater errors compared to the academic example, owing to their increased complexity. The cardiovascular case exhibits an error below 5$\%$ and the industrial case below 15$\%$. Despite the reduced accuracy, significant speed-up is achieved in both applications, exceeding $10^4$.}

Lastly, we introduced the web platforms Atlas and Argos (Sec. \ref{sec:argos}), designed to provide access to industrial and biomedical cases to a broader audience. These platforms utilize data-driven models for real-time predictions, pointing out the potentiality of data-driven ROMs in facilitating user-friendly access to complex modeling tasks by not expert people.

Overall, this chapter shows the applicability of ROMs approaches to different real-world scenarios as well as their capability to be integrated in new technological pipelines. 

\section*{Ackwnoledgments}
We acknowledge the support provided by the European Research Council Executive Agency by the Consolidator Grant project AROMA-CFD ``Advanced Reduced Order Methods with Applications in Computational Fluid Dynamics'' - GA 681447, H2020-ERC CoG 2015 AROMA-CFD, PI G. Rozza, and INdAM-GNCS 2019–2020 projects. P.C.A. and G.R. have also received support by the consortium iNEST (Interconnected North-East Innovation Ecosystem), Piano Nazionale di Ripresa e Resilienza (PNRR) – Missione 4 Componente 2, Investimento 1.5 – D.D. 1058 23/06/2022, ECS\_00000043, supported by the European Union's NextGenerationEU program.

\bibliography{bib.bib}

\begin{thebibliography}{106}
\providecommand{\natexlab}[1]{#1}
\providecommand{\url}[1]{\texttt{#1}}
\expandafter\ifx\csname urlstyle\endcsname\relax
  \providecommand{\doi}[1]{doi: #1}\else
  \providecommand{\doi}{doi: \begingroup \urlstyle{rm}\Url}\fi

\bibitem[Africa(2022)]{africa2022lifex}
P.~C. Africa.
\newblock life\textsuperscript{x}: A flexible, high performance library for the
  numerical solution of complex finite element problems.
\newblock \emph{SoftwareX}, 20:\penalty0 101252, 2022.
\newblock \doi{10.1016/j.softx.2022.101252}.

\bibitem[Africa et~al.(2023{\natexlab{a}})Africa, Perotto, and
  de~Falco]{africa2023quadtree}
P.~C. Africa, S.~Perotto, and C.~de~Falco.
\newblock Scalable recovery-based adaptation on cartesian quadtree meshes for
  advection-diffusion-reaction problems.
\newblock \emph{Advances in Computational Science and Engineering}, 1\penalty0
  (4):\penalty0 443--473, 2023{\natexlab{a}}.
\newblock \doi{10.3934/acse.2023018}.

\bibitem[Africa et~al.(2023{\natexlab{b}})Africa, Piersanti, Fedele, Dede', and
  Quarteroni]{africa2023lifex_fibers}
P.~C. Africa, R.~Piersanti, M.~Fedele, L.~Dede', and A.~Quarteroni.
\newblock life\textsuperscript{x}-fiber: an open tool for myofibers generation
  in cardiac computational models.
\newblock \emph{BMC Bioinformatics}, 24\penalty0 (1):\penalty0 143,
  2023{\natexlab{b}}.
\newblock \doi{10.1186/s12859-023-05260-w}.

\bibitem[Africa et~al.(2023{\natexlab{c}})Africa, Piersanti, Regazzoni,
  Bucelli, Salvador, Fedele, Pagani, Dede', and Quarteroni]{africa2023lifex_ep}
P.~C. Africa, R.~Piersanti, F.~Regazzoni, M.~Bucelli, M.~Salvador, M.~Fedele,
  S.~Pagani, L.~Dede', and A.~Quarteroni.
\newblock life\textsuperscript{x}-ep: a robust and efficient software for
  cardiac electrophysiology simulations.
\newblock \emph{BMC Bioinformatics}, 24\penalty0 (1):\penalty0 389,
  2023{\natexlab{c}}.
\newblock \doi{10.1186/s12859-023-05513-8}.

\bibitem[Africa et~al.(2023{\natexlab{d}})Africa, Salvador, Gervasio, Dede',
  and Quarteroni]{africa2023epmf}
P.~C. Africa, M.~Salvador, P.~Gervasio, L.~Dede', and A.~Quarteroni.
\newblock A matrix--free high--order solver for the numerical solution of
  cardiac electrophysiology.
\newblock \emph{Journal of Computational Physics}, 478:\penalty0 111984,
  2023{\natexlab{d}}.
\newblock \doi{10.1016/j.jcp.2023.111984}.

\bibitem[Africa et~al.(2024)Africa, Fumagalli, Bucelli, Zingaro, Fedele,
  Quarteroni, et~al.]{africa2024lifex_cfd}
P.~C. Africa, I.~Fumagalli, M.~Bucelli, A.~Zingaro, M.~Fedele, A.~Quarteroni,
  et~al.
\newblock life\textsuperscript{x}-cfd: An open-source computational fluid
  dynamics solver for cardiovascular applications.
\newblock \emph{Computer Physics Communications}, 296:\penalty0 109039, 2024.
\newblock \doi{10.1016/j.cpc.2023.109039}.

\bibitem[Amoiralis and Nikolos(2008)]{amoiralis2008freeform}
E.~I. Amoiralis and I.~K. Nikolos.
\newblock {Freeform deformation versus B-spline representation in inverse
  airfoil design}.
\newblock \emph{Journal of Computing and Information Science in Engineering},
  8\penalty0 (2):\penalty0 024001, 04 2008.

\bibitem[Amsallem and Farhat(2012)]{amsallem2012stabilization}
D.~Amsallem and C.~Farhat.
\newblock {Stabilization of projection-based reduced-order models}.
\newblock \emph{International Journal for Numerical Methods in Engineering},
  91\penalty0 (4):\penalty0 358--377, 2012.

\bibitem[Balajewicz and Dowell(2012)]{balajewicz2012stabilization}
M.~Balajewicz and E.~H. Dowell.
\newblock {Stabilization of projection-based reduced order models of the
  Navier--Stokes}.
\newblock \emph{Nonlinear Dynamics}, 70:\penalty0 1619--1632, 2012.

\bibitem[Balzotti et~al.(2022)Balzotti, Siena, Girfoglio, Quaini, and
  Rozza]{balzottidata2022}
C.~Balzotti, P.~Siena, M.~Girfoglio, A.~Quaini, and G.~Rozza.
\newblock {A data-driven reduced order method for parametric optimal blood flow
  control: Application to coronary bypass graft}.
\newblock \emph{Communications in Optimization Theory}, 2022(26):\penalty0
  1--19, 2022.

\bibitem[Balzotti et~al.(2023)Balzotti, Siena, Girfoglio, Stabile,
  Due{\~n}as-Pamplona, Sierra-Pallares, Amat-Santos, and
  Rozza]{balzotti2023reduced}
C.~Balzotti, P.~Siena, M.~Girfoglio, G.~Stabile, J.~Due{\~n}as-Pamplona,
  J.~Sierra-Pallares, I.~Amat-Santos, and G.~Rozza.
\newblock {A Reduced Order Model formulation for left atrium flow: an Atrial
  Fibrillation case}.
\newblock \emph{arXiv preprint arXiv:2309.10601}, 2023.

\bibitem[Bang-Jensen et~al.(2004)Bang-Jensen, Gutin, and Yeo]{bang2004greedy}
J.~Bang-Jensen, G.~Gutin, and A.~Yeo.
\newblock {When the greedy algorithm fails}.
\newblock \emph{Discrete optimization}, 1\penalty0 (2):\penalty0 121--127,
  2004.

\bibitem[Benner et~al.(2020{\natexlab{a}})Benner, Schilders, Grivet-Talocia,
  Quarteroni, Rozza, and Miguel~Silveira]{benner2020model}
P.~Benner, W.~Schilders, S.~Grivet-Talocia, A.~Quarteroni, G.~Rozza, and
  L.~Miguel~Silveira.
\newblock \emph{{Model Order Reduction: Volume 3 Applications}}, volume~3.
\newblock De Gruyter, 2020{\natexlab{a}}.

\bibitem[Benner et~al.(2020{\natexlab{b}})Benner, Schilders, Grivet-Talocia,
  Quarteroni, Rozza, and Miguel~Silveira]{benner2020model2}
P.~Benner, W.~Schilders, S.~Grivet-Talocia, A.~Quarteroni, G.~Rozza, and
  L.~Miguel~Silveira.
\newblock \emph{{Model Order Reduction: Volume 2: Snapshot-Based Methods and
  Algorithms}}.
\newblock De Gruyter, 2020{\natexlab{b}}.

\bibitem[Benner et~al.(2021)Benner, Grivet-Talocia, Quarteroni, Rozza,
  Schilders, and Silveira]{degruyter1}
P.~Benner, S.~Grivet-Talocia, A.~Quarteroni, G.~Rozza, W.~Schilders, and L.~M.
  Silveira.
\newblock \emph{{Model Order Reduction: Volume 1: System- and Data-Driven
  Methods and Algorithms}}.
\newblock De Gruyter, 2021.

\bibitem[Boukraichi et~al.(2023)Boukraichi, Razaaly, Akkari, Casenave, and
  Ryckelynck]{boukraichi2023parametrized}
H.~Boukraichi, N.~Razaaly, N.~Akkari, F.~Casenave, and D.~Ryckelynck.
\newblock {Parametrized non intrusive space-time approximation for explicit
  dynamic fem applications}.
\newblock \emph{ESAIM: Proceedings and Surveys}, 73:\penalty0 68--88, 2023.

\bibitem[Brujic et~al.(2002)Brujic, Ristic, and
  Ainsworth]{brujic2002measurement}
D.~Brujic, M.~Ristic, and I.~Ainsworth.
\newblock {Measurement-based modification of NURBS surfaces}.
\newblock \emph{Computer-Aided Design}, 34\penalty0 (3):\penalty0 173--183,
  2002.

\bibitem[Bucelli et~al.(2023)Bucelli, Zingaro, Africa, Fumagalli, Dede', and
  Quarteroni]{bucelli2023}
M.~Bucelli, A.~Zingaro, P.~C. Africa, I.~Fumagalli, L.~Dede', and
  A.~Quarteroni.
\newblock A mathematical model that integrates cardiac electrophysiology,
  mechanics and fluid dynamics: Application to the human left heart.
\newblock \emph{International Journal for Numerical Methods in Biomedical
  Engineering}, 39\penalty0 (3):\penalty0 e3678, 2023.
\newblock \doi{10.1002/cnm.3678}.

\bibitem[Buhmann and Dyn(1993)]{buhmann1993spectral}
M.~Buhmann and N.~Dyn.
\newblock {Spectral convergence of multiquadric interpolation}.
\newblock \emph{Proceedings of the Edinburgh Mathematical Society}, 36\penalty0
  (2):\penalty0 319--333, 1993.

\bibitem[Buhmann(2000)]{buhmann2000radial}
M.~D. Buhmann.
\newblock {Radial basis functions}.
\newblock \emph{Acta numerica}, 9:\penalty0 1--38, 2000.

\bibitem[Calin(2020)]{calin2020deep}
O.~Calin.
\newblock \emph{{Deep learning architectures}}.
\newblock Springer, 2020.

\bibitem[Chen et~al.(2021)Chen, Wang, Hesthaven, and Zhang]{chen2021physics}
W.~Chen, Q.~Wang, J.~S. Hesthaven, and C.~Zhang.
\newblock {Physics-informed machine learning for reduced-order modeling of
  nonlinear problems}.
\newblock \emph{Journal of Computational Physics}, 446:\penalty0 110666, 2021.

\bibitem[Chen et~al.(2012)Chen, Hesthaven, Maday, Rodr{\'\i}guez, and
  Zhu]{chen2012certified}
Y.~Chen, J.~S. Hesthaven, Y.~Maday, J.~Rodr{\'\i}guez, and X.~Zhu.
\newblock {Certified reduced basis method for electromagnetic scattering and
  radar cross section estimation}.
\newblock \emph{Computer Methods in Applied Mechanics and Engineering},
  233:\penalty0 92--108, 2012.

\bibitem[Coscia et~al.(2023)Coscia, Ivagnes, Demo, and
  Rozza]{coscia2023physics}
D.~Coscia, A.~Ivagnes, N.~Demo, and G.~Rozza.
\newblock {Physics-informed neural networks for advanced modeling}.
\newblock \emph{Journal of Open Source Software}, 8\penalty0 (87):\penalty0
  5352, 2023.

\bibitem[Crowe et~al.(2011)Crowe, Schwarzkopf, Sommerfeld, and Tsuji]{clayton}
C.~T. Crowe, J.~D. Schwarzkopf, M.~Sommerfeld, and Y.~Tsuji.
\newblock \emph{Multiphase {{F}}lows with {{D}}roplets and {{P}}articles},
  volume~10.
\newblock CRC Press, 2011.

\bibitem[Cunningham and Ghahramani(2015)]{cunningham2015linear}
J.~P. Cunningham and Z.~Ghahramani.
\newblock {Linear dimensionality reduction: Survey, insights, and
  generalizations}.
\newblock \emph{The Journal of Machine Learning Research}, 16\penalty0
  (1):\penalty0 2859--2900, 2015.

\bibitem[Cybenko(1989)]{cybenko1989approximation}
G.~Cybenko.
\newblock {Approximation by superpositions of a sigmoidal function}.
\newblock \emph{Mathematics of Control, Signals and Systems}, 2\penalty0
  (4):\penalty0 303--314, 1989.

\bibitem[Daniel et~al.(2020)Daniel, Casenave, Akkari, and
  Ryckelynck]{daniel2020model}
T.~Daniel, F.~Casenave, N.~Akkari, and D.~Ryckelynck.
\newblock {Model order reduction assisted by deep neural networks (ROM-net)}.
\newblock \emph{Advanced Modeling and Simulation in Engineering Sciences},
  7:\penalty0 1--27, 2020.

\bibitem[Dar et~al.(2023)Dar, Baiges, and Codina]{dar2023artificial}
Z.~Dar, J.~Baiges, and R.~Codina.
\newblock {Artificial neural network based correction for reduced order models
  in computational fluid mechanics}.
\newblock \emph{Computer Methods in Applied Mechanics and Engineering},
  415:\penalty0 116232, 2023.

\bibitem[Eckart and Young(1936)]{eckart1936approximation}
C.~Eckart and G.~Young.
\newblock {The approximation of one matrix by another of lower rank}.
\newblock \emph{Psychometrika}, 1\penalty0 (3):\penalty0 211--218, 1936.

\bibitem[Fedele et~al.(2023)Fedele, Piersanti, Regazzoni, Salvador, Africa,
  Bucelli, Zingaro, Dede', and Quarteroni]{fedele2023}
M.~Fedele, R.~Piersanti, F.~Regazzoni, M.~Salvador, P.~C. Africa, M.~Bucelli,
  A.~Zingaro, L.~Dede', and A.~Quarteroni.
\newblock A comprehensive and biophysically detailed computational model of the
  whole human heart electromechanics.
\newblock \emph{Computer Methods in Applied Mechanics and Engineering},
  410:\penalty0 115983, 2023.
\newblock \doi{10.1016/j.cma.2023.115983}.

\bibitem[Feng et~al.(2024)Feng, Chellappa, and Benner]{feng2024posteriori}
L.~Feng, S.~Chellappa, and P.~Benner.
\newblock {A posteriori error estimation for model order reduction of
  parametric systems}.
\newblock \emph{Advanced Modeling and Simulation in Engineering Sciences},
  11\penalty0 (1):\penalty0 5, 2024.

\bibitem[Fernandes et~al.(2018)Fernandes, Semyonov, Ferrás, and
  Nóbrega]{fernandes}
C.~Fernandes, D.~Semyonov, L.~L. Ferrás, and J.~M. Nóbrega.
\newblock Validation of the {{CFD-DPM}} solver {{DPMF}}oam in {{O}}pen{{FOAM}}
  through analytical, numerical and experimental comparisons.
\newblock \emph{Granular Matter}, 20\penalty0 (4):\penalty0 1--18, 2018.

\bibitem[Fine(2006)]{fine2006feedforward}
T.~L. Fine.
\newblock \emph{{Feedforward neural network methodology}}.
\newblock Springer Science \& Business Media, 2006.

\bibitem[Forti and Rozza(2014)]{forti2014efficient}
D.~Forti and G.~Rozza.
\newblock {Efficient geometrical parametrisation techniques of interfaces for
  reduced-order modelling: application to fluid--structure interaction coupling
  problems}.
\newblock \emph{International Journal of Computational Fluid Dynamics},
  28\penalty0 (3-4):\penalty0 158--169, 2014.

\bibitem[Franke(1982)]{franke1982scattered}
R.~Franke.
\newblock {Scattered data interpolation: tests of some methods}.
\newblock \emph{Mathematics of Computation}, 38\penalty0 (157):\penalty0
  181--200, 1982.

\bibitem[Fresca and Manzoni(2022)]{fresca2022pod}
S.~Fresca and A.~Manzoni.
\newblock {POD-DL-ROM: Enhancing deep learning-based reduced order models for
  nonlinear parametrized PDEs by proper orthogonal decomposition}.
\newblock \emph{Computer Methods in Applied Mechanics and Engineering},
  388:\penalty0 114181, 2022.

\bibitem[Fu et~al.(2021)Fu, Xiao, Navon, and Wang]{fu2021data}
R.~Fu, D.~Xiao, I.~Navon, and C.~Wang.
\newblock {A data driven reduced order model of fluid flow by auto-encoder and
  self-attention deep learning methods}.
\newblock \emph{arXiv preprint arXiv:2109.02126}, 2021.

\bibitem[Girfoglio et~al.(2021{\natexlab{a}})Girfoglio, Quaini, and
  Rozza]{girfoglio2021pod}
M.~Girfoglio, A.~Quaini, and G.~Rozza.
\newblock {A POD-Galerkin reduced order model for a LES filtering approach}.
\newblock \emph{Journal of Computational Physics}, 436:\penalty0 110260,
  2021{\natexlab{a}}.

\bibitem[Girfoglio et~al.(2021{\natexlab{b}})Girfoglio, Quaini, and
  Rozza]{girfoglio2021pressure}
M.~Girfoglio, A.~Quaini, and G.~Rozza.
\newblock {Pressure stabilization strategies for a LES filtering reduced order
  model}.
\newblock \emph{Fluids}, 6\penalty0 (9):\penalty0 302, 2021{\natexlab{b}}.

\bibitem[Goldschmidt(2001)]{goldschmidt2001hydrodynamic}
M.~Goldschmidt.
\newblock \emph{Hydrodynamic modelling of fluidised bed spray granulation}.
\newblock PhD thesis, University of Twente, 2001.

\bibitem[Golshan et~al.(2020)Golshan, Sotudeh-Gharebagh, Zarghami, Mostoufi,
  and Blais]{review}
S.~Golshan, R.~Sotudeh-Gharebagh, R.~Zarghami, N.~Mostoufi, and B.~Blais.
\newblock Review and {{I}}mplementation of {{CFD-DEM}} {{A}}pplied to
  {{C}}hemical {{P}}rocess {{S}}ystems.
\newblock \emph{Chemical Engineering Science}, 2020.

\bibitem[Gonzalez and Balajewicz(2018)]{gonzalez2018deep}
F.~J. Gonzalez and M.~Balajewicz.
\newblock {Deep convolutional recurrent autoencoders for learning
  low-dimensional feature dynamics of fluid systems}.
\newblock \emph{arXiv preprint arXiv:1808.01346}, 2018.

\bibitem[Goodfellow et~al.(2016)Goodfellow, Bengio, and
  Courville]{goodfellow2016deep}
I.~Goodfellow, Y.~Bengio, and A.~Courville.
\newblock \emph{{Deep learning}}.
\newblock MIT press, 2016.

\bibitem[Gunzburger et~al.(2007)Gunzburger, Peterson, and
  Shadid]{gunzburger2007reduced}
M.~D. Gunzburger, J.~S. Peterson, and J.~N. Shadid.
\newblock {Reduced-order modeling of time-dependent PDEs with multiple
  parameters in the boundary data}.
\newblock \emph{Computer Methods in Applied Mechanics and Engineering},
  196\penalty0 (4-6):\penalty0 1030--1047, 2007.

\bibitem[Hajisharifi et~al.(2023)Hajisharifi, Roman{\`o}, Girfoglio, Beccari,
  Bonanni, and Rozza]{hajisharifi2023non}
A.~Hajisharifi, F.~Roman{\`o}, M.~Girfoglio, A.~Beccari, D.~Bonanni, and
  G.~Rozza.
\newblock {A non-intrusive data-driven reduced order model for parametrized
  CFD-DEM numerical simulations}.
\newblock \emph{Journal of Computational Physics}, 491:\penalty0 112355, 2023.

\bibitem[Hardy(1971)]{hardy1971multiquadric}
R.~L. Hardy.
\newblock {Multiquadric equations of topography and other irregular surfaces}.
\newblock \emph{Journal of Geophysical Research}, 76\penalty0 (8):\penalty0
  1905--1915, 1971.

\bibitem[Hesthaven and Ubbiali(2018)]{hesthaven2018non}
J.~S. Hesthaven and S.~Ubbiali.
\newblock {Non-intrusive reduced order modeling of nonlinear problems using
  neural networks}.
\newblock \emph{Journal of Computational Physics}, 363:\penalty0 55--78, 2018.

\bibitem[Hesthaven et~al.(2016)Hesthaven, Rozza, Stamm,
  et~al.]{hesthaven2016certified}
J.~S. Hesthaven, G.~Rozza, B.~Stamm, et~al.
\newblock \emph{{Certified reduced basis methods for parametrized partial
  differential equations}}, volume 590.
\newblock Springer, 2016.

\bibitem[Hesthaven et~al.(2022)Hesthaven, Pagliantini, and
  Rozza]{hesthaven2022reduced}
J.~S. Hesthaven, C.~Pagliantini, and G.~Rozza.
\newblock {Reduced basis methods for time-dependent problems}.
\newblock \emph{Acta Numerica}, 31:\penalty0 265--345, 2022.

\bibitem[Hijazi et~al.(2020)Hijazi, Stabile, Mola, and Rozza]{hijazi2020data}
S.~Hijazi, G.~Stabile, A.~Mola, and G.~Rozza.
\newblock {Data-driven POD-Galerkin reduced order model for turbulent flows}.
\newblock \emph{Journal of Computational Physics}, 416:\penalty0 109513, 2020.

\bibitem[Iollo et~al.(2000)Iollo, Lanteri, and
  D{\'e}sid{\'e}ri]{iollo2000stability}
A.~Iollo, S.~Lanteri, and J.-A. D{\'e}sid{\'e}ri.
\newblock {Stability properties of POD--Galerkin approximations for the
  compressible Navier--Stokes equations}.
\newblock \emph{Theoretical and Computational Fluid Dynamics}, 13\penalty0
  (6):\penalty0 377--396, 2000.

\bibitem[Ishida et~al.(2001)Ishida, Sakuma, Cruz, Shimono, Tokui, Yada, Takeda,
  and Higgins]{ishida2001mr}
N.~Ishida, H.~Sakuma, B.~P. Cruz, T.~Shimono, T.~Tokui, I.~Yada, K.~Takeda, and
  C.~B. Higgins.
\newblock {MR flow measurement in the internal mammary artery--to--coronary
  artery bypass graft: comparison with graft stenosis at radiographic
  angiography}.
\newblock \emph{Radiology}, 220\penalty0 (2):\penalty0 441--447, 2001.

\bibitem[Ito and Ravindran(1998)]{ito1998reduced}
K.~Ito and S.~S. Ravindran.
\newblock {A reduced-order method for simulation and control of fluid flows}.
\newblock \emph{Journal of Computational Physics}, 143\penalty0 (2):\penalty0
  403--425, 1998.

\bibitem[Karatzas et~al.(2020)Karatzas, Ballarin, and
  Rozza]{karatzas2020projection}
E.~N. Karatzas, F.~Ballarin, and G.~Rozza.
\newblock {Projection-based reduced order models for a cut finite element
  method in parametrized domains}.
\newblock \emph{Computers \& Mathematics with Applications}, 79\penalty0
  (3):\penalty0 833--851, 2020.

\bibitem[Kashima(2016)]{kashima2016nonlinear}
K.~Kashima.
\newblock {Nonlinear model reduction by deep autoencoder of noise response
  data}.
\newblock In \emph{55th conference on decision and control}, pages 5750--5755.
  IEEE, 2016.

\bibitem[Keegan et~al.(2004)Keegan, Gatehouse, Yang, and
  Firmin]{keegan2004spiral}
J.~Keegan, P.~D. Gatehouse, G.-Z. Yang, and D.~N. Firmin.
\newblock {Spiral phase velocity mapping of left and right coronary artery
  blood flow: Correction for through-plane motion using selective fat-only
  excitation}.
\newblock \emph{Journal of Magnetic Resonance Imaging: An Official Journal of
  the International Society for Magnetic Resonance in Medicine}, 20\penalty0
  (6):\penalty0 953--960, 2004.

\bibitem[Kriesel(2007)]{kriesel2007brief}
D.~Kriesel.
\newblock {A brief introduction to neural networks}, 2007.

\bibitem[Lamousin and Waggenspack(1994)]{lamousin1994nurbs}
H.~J. Lamousin and N.~Waggenspack.
\newblock {NURBS-based free-form deformations}.
\newblock \emph{IEEE Computer Graphics and Applications}, 14\penalty0
  (6):\penalty0 59--65, 1994.

\bibitem[Lassila et~al.(2013)Lassila, Manzoni, Quarteroni, and
  Rozza]{lassila2013generalized}
T.~Lassila, A.~Manzoni, A.~Quarteroni, and G.~Rozza.
\newblock {Generalized reduced basis methods and n-width estimates for the
  approximation of the solution manifold of parametric PDEs}.
\newblock In \emph{Analysis and numerics of partial differential equations},
  pages 307--329. Springer, 2013.

\bibitem[Lee and Carlberg(2020)]{lee2020model}
K.~Lee and K.~T. Carlberg.
\newblock {Model reduction of dynamical systems on nonlinear manifolds using
  deep convolutional autoencoders}.
\newblock \emph{Journal of Computational Physics}, 404:\penalty0 108973, 2020.

\bibitem[Lindner et~al.(2017)Lindner, Mehl, and Uekermann]{lindner2017radial}
F.~Lindner, M.~Mehl, and B.~Uekermann.
\newblock {Radial basis function interpolation for black-box multi-physics
  simulations}.
\newblock In \emph{VII International Conference on Computational Methods for
  Coupled Problems in Science and Engineering}. CIMNE, 2017.

\bibitem[Milano and Koumoutsakos(2002)]{milano2002neural}
M.~Milano and P.~Koumoutsakos.
\newblock {Neural network modeling for near wall turbulent flow}.
\newblock \emph{Journal of Computational Physics}, 182\penalty0 (1):\penalty0
  1--26, 2002.

\bibitem[Moliner et~al.(2018)Moliner, Marchelli, Spanachi, Martinez-Felipe,
  Bosio, and Arato]{molin}
C.~Moliner, F.~Marchelli, N.~Spanachi, A.~Martinez-Felipe, B.~Bosio, and
  E.~Arato.
\newblock {{CFD}} simulation of a spouted bed: Comparison between the
  {{D}}iscrete {{E}}lement {{M}}ethod ({{DEM}}) and the {{T}}wo {{F}}luid
  {{M}}odel ({{TFM}}).
\newblock \emph{Chemical Engineering Journal}, 2018.

\bibitem[Monk et~al.(2003)]{monk2003finite}
P.~Monk et~al.
\newblock \emph{{Finite element methods for Maxwell's equations}}.
\newblock Oxford University Press, 2003.

\bibitem[Nguyen(2024)]{nguyen2024model}
N.~C. Nguyen.
\newblock {Model reduction techniques for parametrized nonlinear partial
  differential equations}.
\newblock \emph{Advances in Applied Mechanics (AAMS)}, 58, 2024.

\bibitem[Nurtaj~Hossain and Ghosh(2020)]{nurtaj2020adaptive}
M.~Nurtaj~Hossain and D.~Ghosh.
\newblock {Adaptive reduced order modeling for nonlinear dynamical systems
  through a new a posteriori error estimator: Application to uncertainty
  quantification}.
\newblock \emph{International Journal for Numerical Methods in Engineering},
  121\penalty0 (15):\penalty0 3417--3441, 2020.

\bibitem[Papapicco et~al.(2022)Papapicco, Demo, Girfoglio, Stabile, and
  Rozza]{papapicco2021neural}
D.~Papapicco, N.~Demo, M.~Girfoglio, G.~Stabile, and G.~Rozza.
\newblock {The Neural Network shifted-proper orthogonal decomposition: A
  machine learning approach for non-linear reduction of hyperbolic equations}.
\newblock \emph{Computer Methods in Applied Mechanics and Engineering},
  392:\penalty0 114687, 2022.

\bibitem[Pichi and Rozza(2019)]{pichi2019reduced}
F.~Pichi and G.~Rozza.
\newblock {Reduced basis approaches for parametrized bifurcation problems held
  by non-linear Von K{\'a}rm{\'a}n equations}.
\newblock \emph{Journal of Scientific Computing}, 81:\penalty0 112--135, 2019.

\bibitem[Pichi et~al.(2020)Pichi, Quaini, and Rozza]{pichi2020reduced}
F.~Pichi, A.~Quaini, and G.~Rozza.
\newblock {A reduced order modeling technique to study bifurcating phenomena:
  application to the Gross--Pitaevskii equation}.
\newblock \emph{SIAM Journal on Scientific Computing}, 42\penalty0
  (5):\penalty0 B1115--B1135, 2020.

\bibitem[Pichi et~al.(2022)Pichi, Strazzullo, Ballarin, and
  Rozza]{pichi2022driving}
F.~Pichi, M.~Strazzullo, F.~Ballarin, and G.~Rozza.
\newblock {Driving bifurcating parametrized nonlinear PDEs by optimal control
  strategies: application to Navier--Stokes equations with model order
  reduction}.
\newblock \emph{ESAIM: Mathematical Modelling and Numerical Analysis},
  56\penalty0 (4):\penalty0 1361--1400, 2022.

\bibitem[Pichi et~al.(2023)Pichi, Ballarin, Rozza, and
  Hesthaven]{pichi2023artificial}
F.~Pichi, F.~Ballarin, G.~Rozza, and J.~S. Hesthaven.
\newblock {An artificial neural network approach to bifurcating phenomena in
  computational fluid dynamics}.
\newblock \emph{Computers \& Fluids}, 254:\penalty0 105813, 2023.

\bibitem[Piersanti et~al.(2021)Piersanti, Africa, Fedele, Vergara, Dede',
  Corno, and Quarteroni]{piersanti2021}
R.~Piersanti, P.~C. Africa, M.~Fedele, C.~Vergara, L.~Dede', A.~F. Corno, and
  A.~Quarteroni.
\newblock Modeling cardiac muscle fibers in ventricular and atrial
  electrophysiology simulations.
\newblock \emph{Computer Methods in Applied Mechanics and Engineering},
  373:\penalty0 113468, 2021.
\newblock \doi{10.1016/j.cma.2020.113468}.

\bibitem[Pinkus(2012)]{pinkus2012n}
A.~Pinkus.
\newblock \emph{{N-widths in Approximation Theory}}, volume~7.
\newblock Springer Science \& Business Media, 2012.

\bibitem[Porsching(1985)]{porsching1985estimation}
T.~A. Porsching.
\newblock {Estimation of the error in the reduced basis method solution of
  nonlinear equations}.
\newblock \emph{Mathematics of Computation}, 45\penalty0 (172):\penalty0
  487--496, 1985.

\bibitem[Quarteroni and Valli(2008)]{quarteroni2008numerical}
A.~Quarteroni and A.~Valli.
\newblock \emph{{Numerical approximation of partial differential equations}},
  volume~23.
\newblock Springer Science \& Business Media, 2008.

\bibitem[Quarteroni et~al.(2011)Quarteroni, Rozza, and
  Manzoni]{quarteroni2011certified}
A.~Quarteroni, G.~Rozza, and A.~Manzoni.
\newblock {Certified reduced basis approximation for parametrized partial
  differential equations and applications}.
\newblock \emph{Journal of Mathematics in Industry}, 1:\penalty0 1--49, 2011.

\bibitem[Quarteroni et~al.(2014)Quarteroni, Rozza,
  et~al.]{quarteroni2014reduced}
A.~Quarteroni, G.~Rozza, et~al.
\newblock \emph{{Reduced order methods for modeling and computational
  reduction}}, volume~9.
\newblock Springer, 2014.

\bibitem[Quarteroni et~al.(2015)Quarteroni, Manzoni, and
  Negri]{quarteroni2015reduced}
A.~Quarteroni, A.~Manzoni, and F.~Negri.
\newblock \emph{{Reduced basis methods for partial differential equations: an
  introduction}}, volume~92.
\newblock Springer, 2015.

\bibitem[Regazzoni et~al.(2022)Regazzoni, Salvador, Africa, Fedele, Dede', and
  Quarteroni]{regazzoni2022}
F.~Regazzoni, M.~Salvador, P.~C. Africa, M.~Fedele, L.~Dede', and
  A.~Quarteroni.
\newblock A cardiac electromechanics model coupled with a lumped-parameter
  model for closed-loop blood circulation.
\newblock \emph{Journal of Computational Physics}, 457:\penalty0 111083, 2022.
\newblock \doi{10.1016/j.jcp.2022.111083}.

\bibitem[Rojas and Rojas(1996)]{rojas1996backpropagation}
R.~Rojas and R.~Rojas.
\newblock {The backpropagation algorithm}.
\newblock \emph{Neural Networks: a Systematic Introduction}, pages 149--182,
  1996.

\bibitem[Rosenblatt(1958)]{rosenblatt1958perceptron}
F.~Rosenblatt.
\newblock {The perceptron: a probabilistic model for information storage and
  organization in the brain.}
\newblock \emph{Psychological Review}, 65\penalty0 (6):\penalty0 386, 1958.

\bibitem[Rowley et~al.(2004)Rowley, Colonius, and Murray]{rowley2004model}
C.~W. Rowley, T.~Colonius, and R.~M. Murray.
\newblock {Model reduction for compressible flows using POD and Galerkin
  projection}.
\newblock \emph{Physica D: Nonlinear Phenomena}, 189\penalty0 (1-2):\penalty0
  115--129, 2004.

\bibitem[Rozza et~al.(2008)Rozza, Huynh, and Patera]{rozza2008reduced}
G.~Rozza, D.~B.~P. Huynh, and A.~T. Patera.
\newblock {Reduced basis approximation and a posteriori error estimation for
  affinely parametrized elliptic coercive partial differential equations}.
\newblock \emph{Archives of Computational Methods in Engineering}, 15\penalty0
  (3):\penalty0 229--275, 2008.

\bibitem[Rozza et~al.(2020)Rozza, Hess, Stabile, Tezzele, and
  Ballarin]{rozza20201}
G.~Rozza, M.~Hess, G.~Stabile, M.~Tezzele, and F.~Ballarin.
\newblock {1 Basic ideas and tools for projection-based model reduction of
  parametric partial differential equations}.
\newblock In \emph{Snapshot-based methods and algorithms}, pages 1--47. De
  Gruyter, 2020.

\bibitem[Rozza et~al.(2022)Rozza, Stabile, and Ballarin]{aromabook}
G.~Rozza, G.~Stabile, and F.~Ballarin.
\newblock \emph{{Advanced Reduced Order Methods and Applications in
  Computational Fluid Dynamics}}.
\newblock SIAM, Philadelphia, PA, 2022.

\bibitem[Rumelhart et~al.(1986)Rumelhart, Hinton, and
  Williams]{rumelhart1986learning}
D.~E. Rumelhart, G.~E. Hinton, and R.~J. Williams.
\newblock {Learning representations by back-propagating errors}.
\newblock \emph{nature}, 323\penalty0 (6088):\penalty0 533--536, 1986.

\bibitem[Salvador et~al.(2021)Salvador, Fedele, Africa, Sung, Dede', Prakosa,
  Trayanova, Chrispin, and Quarteroni]{salvador2021}
M.~Salvador, M.~Fedele, P.~C. Africa, E.~Sung, L.~Dede', A.~Prakosa,
  N.~Trayanova, J.~Chrispin, and A.~Quarteroni.
\newblock Electromechanical modeling of human ventricles with ischemic
  cardiomyopathy: numerical simulations in sinus rhythm and under arrhythmia.
\newblock \emph{Computers in Biology and Medicine}, 136:\penalty0 104674, 2021.
\newblock \doi{10.1016/j.compbiomed.2021.104674}.

\bibitem[Shah et~al.(2021)Shah, Girfoglio, Quintela, Rozza, Lengomin, Ballarin,
  and Barral]{shah2021finite}
N.~V. Shah, M.~Girfoglio, P.~Quintela, G.~Rozza, A.~Lengomin, F.~Ballarin, and
  P.~Barral.
\newblock {Finite element based model order reduction for parametrized one-way
  coupled steady state linear thermomechanical problems}.
\newblock \emph{arXiv preprint arXiv:2111.08534}, 2021.

\bibitem[Sharma et~al.(2017)Sharma, Sharma, and Athaiya]{sharma2017activation}
S.~Sharma, S.~Sharma, and A.~Athaiya.
\newblock {Activation functions in neural networks}.
\newblock \emph{Towards Data Sci}, 6\penalty0 (12):\penalty0 310--316, 2017.

\bibitem[Siena et~al.(2023{\natexlab{a}})Siena, Girfoglio, Ballarin, and
  Rozza]{Siena2022}
P.~Siena, M.~Girfoglio, F.~Ballarin, and G.~Rozza.
\newblock {Data-driven reduced order modelling for patient-specific
  hemodynamics of coronary artery bypass grafts with physical and geometrical
  parameters}.
\newblock \emph{Journal of Scientific Computing}, 94\penalty0 (2):\penalty0 38,
  2023{\natexlab{a}}.

\bibitem[Siena et~al.(2023{\natexlab{b}})Siena, Girfoglio, and
  Rozza]{siena2023fast}
P.~Siena, M.~Girfoglio, and G.~Rozza.
\newblock {Fast and accurate numerical simulations for the study of coronary
  artery bypass grafts by artificial neural networks}.
\newblock In \emph{Reduced Order Models for the Biomechanics of Living Organs},
  pages 167--183. Elsevier, 2023{\natexlab{b}}.

\bibitem[Siena et~al.(2023{\natexlab{c}})Siena, Girfoglio, and
  Rozza]{siena2023introduction}
P.~Siena, M.~Girfoglio, and G.~Rozza.
\newblock {An introduction to POD-greedy-Galerkin reduced basis method}.
\newblock In \emph{Reduced Order Models for the Biomechanics of Living Organs},
  pages 127--145. Elsevier, 2023{\natexlab{c}}.

\bibitem[Skala(2016)]{skala2016practical}
V.~Skala.
\newblock {A practical use of radial basis functions interpolation and
  approximation}.
\newblock \emph{Investigacion Operacional}, 37\penalty0 (2):\penalty0 137--144,
  2016.

\bibitem[Stabile and Rozza(2018)]{stabile2018finite}
G.~Stabile and G.~Rozza.
\newblock {Finite volume POD-Galerkin stabilised reduced order methods for the
  parametrised incompressible Navier--Stokes equations}.
\newblock \emph{Computers \& Fluids}, 173:\penalty0 273--284, 2018.

\bibitem[Stabile et~al.(2017)Stabile, Hijazi, Mola, Lorenzi, and
  Rozza]{stabile2017pod}
G.~Stabile, S.~Hijazi, A.~Mola, S.~Lorenzi, and G.~Rozza.
\newblock {POD-Galerkin reduced order methods for CFD using Finite Volume
  Discretisation: vortex shedding around a circular cylinder}.
\newblock \emph{Communications in Applied and Industrial Mathematics},
  8:\penalty0 210--236, 2017.

\bibitem[Stella et~al.(2020)Stella, Vergara, Maines, Catanzariti, Africa,
  Dematt\`{e}, Centonze, Nobile, Del~Greco, and Quarteroni]{stella2020}
S.~Stella, C.~Vergara, M.~Maines, D.~Catanzariti, P.~C. Africa, C.~Dematt\`{e},
  M.~Centonze, F.~Nobile, M.~Del~Greco, and A.~Quarteroni.
\newblock Integration of activation maps of epicardial veins in computational
  cardiac electrophysiology.
\newblock \emph{Computers in Biology and Medicine}, 127:\penalty0 104047, 2020.
\newblock \doi{10.1016/j.compbiomed.2020.104047}.

\bibitem[Tsuji et~al.(1993)Tsuji, Kawaguchi, and Tanaka]{tsu}
Y.~Tsuji, T.~Kawaguchi, and T.~Tanaka.
\newblock Discrete particle simulation of two-dimensional fluidized bed.
\newblock \emph{Powder Technol.}, 77:\penalty0 79--87, 1993.

\bibitem[Vergara et~al.(2022)Vergara, Stella, Maines, Africa, Catanzariti,
  Dematt\`{e}, Centonze, Nobile, Quarteroni, and Del~Greco]{vergara2022}
C.~Vergara, S.~Stella, M.~Maines, P.~C. Africa, D.~Catanzariti, C.~Dematt\`{e},
  M.~Centonze, F.~Nobile, A.~Quarteroni, and M.~Del~Greco.
\newblock Computational electrophysiology of the coronary sinus branches based
  on electro-anatomical mapping for the prediction of the latest activated
  region.
\newblock \emph{Medical \& Biological Engineering \& Computing}, 2022.
\newblock \doi{10.1007/s11517-022-02610-3}.

\bibitem[Verim et~al.(2015)Verim, {\"O}zt{\"u}rk, K{\"u}{\c{c}}{\"u}k, Kara,
  Sa{\u{g}}lam, and Karde{\c{s}}o{\u{g}}lu]{verim2015cross}
S.~Verim, E.~{\"O}zt{\"u}rk, U.~K{\"u}{\c{c}}{\"u}k, K.~Kara, M.~Sa{\u{g}}lam,
  and E.~Karde{\c{s}}o{\u{g}}lu.
\newblock {Cross-sectional area measurement of the coronary arteries using CT
  angiography at the level of the bifurcation: is there a relationship?}
\newblock \emph{Diagnostic and interventional radiology}, 21\penalty0
  (6):\penalty0 454, 2015.

\bibitem[Wang et~al.(2019)Wang, Hesthaven, and Ray]{wang2019non}
Q.~Wang, J.~S. Hesthaven, and D.~Ray.
\newblock {Non-intrusive reduced order modeling of unsteady flows using
  artificial neural networks with application to a combustion problem}.
\newblock \emph{Journal of computational physics}, 384:\penalty0 289--307,
  2019.

\bibitem[Weller et~al.(1998)Weller, Tabor, Jasak, and Fureby]{weller}
H.~G. Weller, G.~Tabor, H.~Jasak, and C.~Fureby.
\newblock A tensorial approach to computational continuum mechanics using
  object-oriented techniques.
\newblock \emph{Computers in Physics}, 12:\penalty0 620–631, 1998.

\bibitem[Xu and Yu(1997)]{Xu}
B.~H. Xu and A.~B. Yu.
\newblock Numerical simulation of the gas-solid flow in a fluidized bed by
  combining discrete particle method with computational fluid dynamics.
\newblock \emph{Chemical Engineering Science}, 52:\penalty0 2785–2809, 1997.

\bibitem[Zainib et~al.(2022)Zainib, Siena, Girfoglio, Hess, Ballarin, and
  Rozza]{zainib2022chapter}
Z.~Zainib, P.~Siena, M.~Girfoglio, M.~W. Hess, F.~Ballarin, and G.~Rozza.
\newblock {Chapter 18: Reduced Order Methods for Hemodynamics Applications}.
\newblock In \emph{Advanced Reduced Order Methods and Applications in
  Computational Fluid Dynamics}, pages 365--377. SIAM, 2022.

\bibitem[Zhu et~al.(2007)Zhu, Zhou, Yang, and Yu]{zhu2007discrete}
H.~P. Zhu, Z.~Y. Zhou, R.~Y. Yang, and A.~B. Yu.
\newblock Discrete particle simulation of particulate systems: theoretical
  developments.
\newblock \emph{Chemical Engineering Science}, 62\penalty0 (13):\penalty0
  3378--3396, 2007.

\bibitem[Zingaro et~al.(2022)Zingaro, Fumagalli, Fedele, Africa, Dede',
  Quarteroni, and Corno]{zingaro2022}
A.~Zingaro, I.~Fumagalli, M.~Fedele, P.~C. Africa, L.~Dede', A.~Quarteroni, and
  A.~F. Corno.
\newblock A geometric multiscale model for the numerical simulation of blood
  flow in the human left heart.
\newblock \emph{Discrete and Continuous Dynamical Systems - S}, 15\penalty0
  (8):\penalty0 2391--2427, 2022.
\newblock \doi{10.3934/dcdss.2022052}.

\end{thebibliography}


%
%
%
%
%
%
%
%
%
\let\clearpage\relax
%
%
%
%
%
%
%
%
%
%
%
%
%
%
%
%
%
\end{document}